\documentclass[opre,nonblindrev]{informs3} 

\usepackage{endnotes}
\let\footnote=\endnote
\let\enotesize=\normalsize
\def\notesname{Endnotes}%
\def\makeenmark{$^{\theenmark}$}
\def\enoteformat{\rightskip0pt\leftskip0pt\parindent=1.75em
  \leavevmode\llap{\theenmark.\enskip}}

\usepackage{natbib}
 \bibpunct[, ]{(}{)}{,}{a}{}{,}%
 \def\bibfont{\small}%
 \def\bibsep{\smallskipamount}%
 \def\bibhang{24pt}%
\TheoremsNumberedThrough     
\ECRepeatTheorems

\EquationsNumberedThrough    

                 \usepackage{epsf}
\usepackage{epstopdf}
\usepackage{amsmath}
\usepackage{amssymb}
\usepackage{amsxtra}
\usepackage{algorithmic}
\usepackage{algorithm}
\usepackage{mathrsfs}
\usepackage{amsfonts}
\usepackage{dsfont}

\usepackage{latexsym}
\usepackage{float}
\usepackage{flushend}
\usepackage{cuted}

\begin{document}

\RUNAUTHOR{Al-Kanj et al.}
\RUNTITLE{The Information-Collecting Vehicle Routing Problem: Stochastic Optimization for Emergency Storm Response}
\TITLE{The Information-Collecting Vehicle Routing Problem: Stochastic Optimization for Emergency Storm Response}
\ARTICLEAUTHORS{
\AUTHOR{Lina Al-Kanj, Warren B. Powell and Belgacem Bouzaiene-Ayari}
\AFF{Operations Research and Financial Engineering Department, Princeton University, NJ, USA \EMAIL{\{lalkanj,powell,belgacem\}@princeton.edu}}
}

\ABSTRACT{
Utilities face the challenge of responding to power outages due to storms and ice damage, but most power grids are not equipped with sensors to pinpoint the precise location of the faults causing the outage.  Instead, utilities have to depend primarily on phone calls (trouble calls) from customers who have lost power to guide the dispatching of utility trucks.  In this paper,  we develop a policy that routes a utility truck to restore outages in the power grid as quickly as possible, using phone calls to create beliefs about outages, but also using utility trucks as a mechanism for collecting additional information.  This means that routing decisions change not only the physical state of the truck (as it moves from one location to another) and the grid (as the truck performs repairs),  but also our belief about the network, creating the first stochastic vehicle routing problem that explicitly models information collection and belief modeling. We address the problem of managing a single utility truck, which we start by formulating as a sequential stochastic optimization model which captures our belief about the state of the grid.  We propose a stochastic lookahead policy, and use Monte Carlo tree search (MCTS) to produce a practical policy that is asymptotically optimal. Simulation results show that the developed policy restores the power grid much faster compared to standard industry heuristics.
}
\SUBJECTCLASS{Transportation: Vehicle Routing, Dynamic Programming: Applications}
\maketitle

\section{Introduction}
Climate change is producing more powerful storms, increasing the frequency and severity of outages in our power grid.  Despite the importance of electricity in every aspect of our lives, we often do not know the location of  the fault causing an outage, which complicates the task of restoring power. Instead, utilities depend heavily on phone calls (known as ``trouble calls") from customers who have lost their power, in addition to information from the tripping of some circuit breakers/protective devices. Complicating the situation is that as few as one percent of customers call when they lose their power, creating tremendous uncertainty in the knowledge of the state of the grid. This in turn complicates the task of dispatching utility trucks to the location of faults (which are unknown to the utility center) to restore the grid as quickly as~possible.

According to the Edison Electric Institute~(\citet{EEI13}), on average $55$\% of power outages in the U.S. are due to weather and it can reach up to 80\% in some years~(\citet{Camp12}). Wind (primarily) and ice frequently bring trees and branches down on power lines creating sporadic outages that quickly spread through the grid due to the limited number of protective devices that are triggered due to a short circuit. The remaining outages are due to equipment faults, vehicle or construction accidents, maintenance or human errors, and animals. It is estimated that $90$\% of customer outage-minutes are due to faults affecting the local distribution systems~(\citet{Camp12}). However, the remaining $10$\% stem from generation and transmission problems, which can cause wider-scale outages affecting larger numbers of~customers.

Whenever a power line is damaged, the closest protective device through which the power passes to the damaged power line disconnects resulting in power outage to all customers served through that protective device. Protective devices, such as fuses or circuit breakers, are usually responsible for monitoring the power flow through all the attached downstream components that include power lines, substations, transformers and other protective devices which can form a major part of the circuit. The protective device disconnects the damaged circuit to avoid overloading of the network components which can severely damage them. Whenever a power outage occurs, the electric utility center (EUC) is typically notified by the phone calls of the customers. However, a fault in one location trips the first upstream protective device but it can create a much wider set of outages, with random phone calls coming from locations far from the fault. Complicating the situation is that only a few customers will call. Consequently, identifying the locations of the faults across the power system is a major dilemma for the EUC which needs to restore the network quickly.

Even if the EUC knows the location of the faults, routing the truck to restore the grid to minimize the number of customers in outage at any point in time is a challenge.  Currently, EUCs use a simple policy to route the utility trucks to restore power such as first-call, first-serve where the customers that call first to report outages are served first by searching around the call location for the fault. This policy might be far from optimal since the actual fault could be a few kilometers away from the location of call. We believe that managing the truck dispatching problem based on a model that properly captures both the physical system as well as the uncertainty offers the potential of significantly reducing the amount of time that customers have lost power.

In~(\citet{ABP14}), we have developed a probability model that estimates the likelihood of faults occurring across the power grid using two sources of information: 1) the prior probability of power line faults depending on the storm pattern, the structure of the grid and the EUC knowledge of the power line conditions and environment, e.g., existing trees that could potentially fall on the power lines and 2) the phone calls of the customers. In this paper,  a stochastic optimization model is designed that makes use of the information provided by the probability model about the likelihood of the location of faults to route the utility efficiently across the power grid with the objective of restoring power as fast as possible. 



This paper makes the following contributions:  1) it provides the first formal model of a stochastic vehicle routing problem that explicitly captures the physical state (the location of the truck and known state of the grid), information state (other information such as the history of the truck and phone calls), and belief state (capturing probabilistic knowledge of potential outages).   2) It proposes a stochastic lookahead policy that is solved using a specialized implementation of Monte Carlo tree search, a method that offers an asymptotically optimal policy, which is also a first for this problem class.    3) The performance of the developed policy  is tested using a simulator that simulates the power  grid of the state of New Jersey by using real data provided by PSE\&G which is the main EUC in New Jersey. Simulations show that a) the lookahead learning policy closely matches the performance of the optimal solution  and b)  significantly outperforms standard industry heuristics when tested on realistic stochastic models. We also demonstrate that the industry heuristic does not make good use of better information, while our method provide information-consistent behavior, with better results as information improves (as one would expect).


The paper is organized as follows. Section~\ref{sec:lit_rev} summarizes the literature of fault identification, utility crew dispatching for restoring power distribution systems and stochastic vehicle routing. Section~\ref{sec:sys_model} gives a mathematical model of the flow of information for a grid as some event such as a storm evolves producing faults and loss of power. Section~\ref{sec:opt_prob} derives the sequential stochastic optimization problem that routes the utility truck across the grid and discusses several dispatch policies. Section~\ref{sec:MCTS} introduces MCTS as the lookahead policy that approximates the optimal one with good computational efficiency. Section~\ref{sec:results} compares the proposed learning policy against an industry-standard escalation policy. 

\section{Literature Review}\label{sec:lit_rev}
We first review the literature from the power community on outage prediction and the dispatching of utility crews, followed by a review of the relevant articles from the mainstream literature on dynamic vehicle routing.

\subsection{Outage Prediction and Utility Crew Dispatching Literature}

Utilities face two problems in responding to storm damage: 1) identifying the location of outages and 2) dispatching the utility crews.

The power grid contains protective devices that detect power flow interruption upon which they shut down power flow to the affected components to avoid further damage. Thus, in case one of the components of the circuit faults, the closest upstream protective device shuts down causing power outage to all the downstream components. Different utilities use different methodologies for management and power restoration during an outage; while most utilities install measurement units on most of the long-distance, high voltage transmission system, distribution systems are not equally automated due to the cost of the equipments. Some utilities,  such as Carolina Power \& Light utility (\citet{Lamp02}), have advanced automated systems that include SCADA (supervisory control and data acquisition) which is a system formed of sensors across the grid operating with coded signals over communication channels so as to provide control of remote equipment. However, the majority of the utilities do not have large scale installed sensors. Instead, utilities depend primarily on phone calls from customers who have lost their power. Complicating the situation is that as few as one percent of customers call when they lose their power, creating tremendous uncertainty in the knowledge of the state of the grid.

Distribution systems have different configurations; the most studied configuration in the literature is the radial distribution system which has a tree structure starting from the substation. In this model, there is a single path for power flow from the substation to the consumers. This enables utilities and researchers to propose escalation algorithms for fault identification based on customer calls and grid topology~(\citet{Scot90,HLC91}). Escalation algorithms gather the set of calls and then from the location of each customer call, the tree is followed  back to the substation until the first common location for all calls is located; this location would be identified as the faulted one. Escalation algorithms suit single fault scenarios but cannot capture the case of multiple fault scenarios. An improved escalation algorithm for  a heat storm is presented in~(\citet{LS99}) where the escalation from the locations of calls depends on the type of the upstream devices. Artificial intelligence techniques that make use of customer calls to identify the fault locations are also investigated to provide better performance than the simple heuristics described above. For example a neural network is presented in~(\citet{LTH94}) but the main limitation is in determining the sample training set. An approach based on fuzzy set theory and tabu search is proposed in~(\citet{CW98}) but the limitation is in the computational complexity that becomes intractable for large networks.  Knowledge-based outage identification that make use of SCADA and automated meter reading to provide the EUC with knowledge about the status of the distribution system on top of customer  calls is proposed in~(\citet{LS02}). But, most of the utilities do not have automated distribution systems and thus, primarily rely on the grid topology, the phone calls of the customers and the experience of the utility personnel to estimate the locations of outages. Moreover, since only a few customers will call, identifying the locations of the faults across the power system is still a major dilemma for the EUC which needs to restore the network as fast as possible.

Managing utility crews to restore outages attracted a modest level of attention in the research literature~(\citet{SLK16,ZDT98,Wh14}). In~(\citet{ZDT98}), a utility truck is routed to the location of a customer that called to report a power outage; then, after restoring the fault, the truck is routed to the next calling customer that has the shortest travel time. In~(\citet{Wh14}), the authors develop a model based on data mining and machine learning techniques to predict outages in the grid using collected data from past six storms  as well as asset information (framing, pole
age, etc.), in addition to environmental information. Then, given the predicted outages, deterministic optimization models are developed to route a given number of utility trucks to perform a predefined number of repair jobs required at each damaged location in order to minimize the grid restoration time. The authors in~(\citet{SLK16}) propose an outage prediction model at the level of areas served by a substation, then a mechanism is proposed to  plan hourly crew staffing levels across different organizations (service centers, local contractors, mutual aid crews) and different crew types in order to minimize the overall grid restoration time but they do not solve the problem of utility crew routing. 

\subsection{Stochastic Vehicle Routing Literature}
The study of stochastic and dynamic vehicle routing has a long history, with survey articles appearing as early as (\citet{StGo1982}).  The problem area has attracted so much attention that there has been a steady flow of reviews and survey articles (\citet{StGo1983, Ps1988, PoJa1995, GePo1998, BeCo2010, PiGe2013, LaMa2014}).  Most of this literature focuses on uncertainty in pickups  or deliveries, although some authors address random travel times (\citet{KeMo2003}).   Algorithmic strategies range from scenario-based lookahead policies (e.g. (\citet{LaHa2002,BeVa2004}), to a wide variety of rolling horizon heuristics that involve solving deterministic approximations modified to handle uncertainty.

Most of the academic literature has focused on solving stochastic lookahead models (\citet{LaHa2002, BeVa2004}).  Since these problems combine uncertainty with the complexity of these difficult integer programs, these lookahead models are themselves quite hard to solve.  Often overlooked is that they are simply rolling horizon heuristics to solve a fully-sequential problem, which is typically not modeled (but is often represented in a simulator). \citet{PoSi2012} provides a general framework for modeling stochastic, dynamic problems in transportation and logistics, identifying four classes of policies.  Focusing specifically on vehicle routing, (\citet{GoOh2013}) provides a proper model of the stochastic vehicle routing problem, and then proposes a broad class of practical roll-out policies.

Our paper addresses the problem of routing a single utility truck that has to find and repair outages over a power grid.  Since this problem is completely new, we continue a long tradition of beginning by addressing the single-vehicle version of the problem (see, for example, (\citet{Ps1980}) and (\citet{ReGe2010})). Our problem includes a unique dimension which has never been addressed in the literature: probabilistic knowledge about the state of the network, and the dimension that utility trucks are not only repairing the network, but as the vehicle progresses, it collects information that is used to update the belief knowledge.  This “state of knowledge” has to be represented along with the physical state (location) of the network, and the history of the truck.

Our solution strategy uses a technique familiar to the computer science community known as Monte Carlo Tree Search (MCTS).  While this is most commonly used for deterministic problems (\citet{BPW12}), it has been extended to stochastic problems and shown to be asymptotically optimal for certain variations of the algorithm (\citet{ACT13}). MCTS has been popular in the computer science community (although primarily for deterministic problems), but has only recently seen applications in transportation and logistics (\citet{MaNe2015, EdGa2016}). However, we are unaware of any prior application of MCTS in the setting of information collection, where the state variable includes a belief state.

\section{Problem Description}~\label{sec:sys_model}

\begin{figure}[t!]
\vspace{-0cm}
\centering
\includegraphics[scale=0.4]{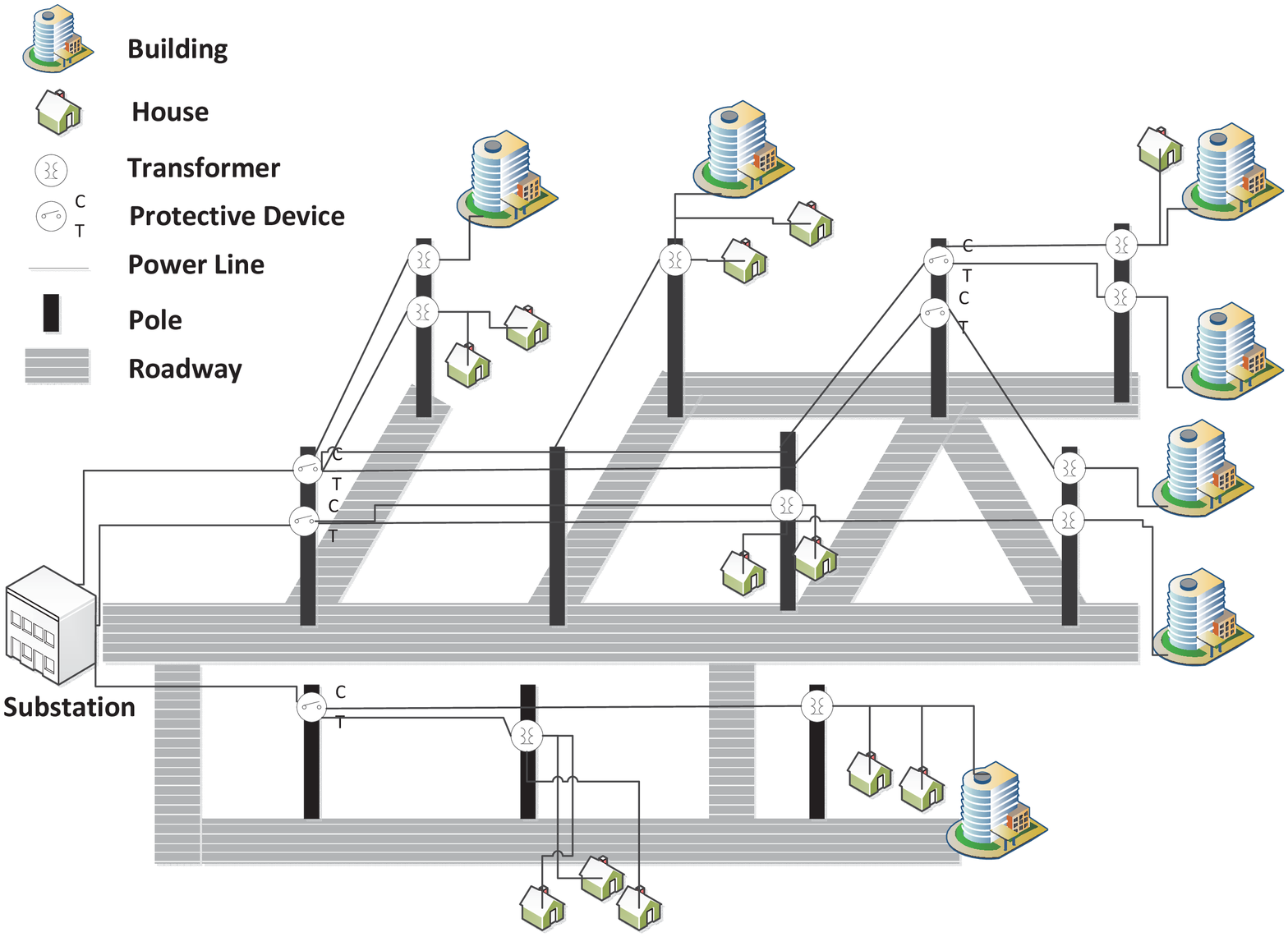}
\vspace{-0.0cm}
\caption{\label{fig:sys_mod} Illustration of a circuit in a distribution grid, with feeders from a single substation.}\vspace{-0.3cm}
\end{figure}

In this section, we describe the parameters and variables that govern the evolution of the problem to be solved over time. First, we assume an overhead power system consisting of substations, poles that carry the protective devices, power lines and transformers through which customers are fed with power as shown in Figure~\ref{fig:sys_mod}. In~(\citet{ABP14}), we have already developed a probability model that estimates the likelihood of fault occurring across the power lines of the grid based on the prior probability of power line fault and the phone calls of the customers.
It is assumed that a fault across a power line includes the fault across its connection ends which could be a transformer or a protective device. This can be justified by the fact that a fault across a power line or one of the components at its end connections results in power outage to the same set of customers in the power system.


Let $\mathcal{U}$ be the set of circuits in the power system, $\mathcal{U}=\{u: u=1,\ldots,U\}$ where each circuit can be graphically represented by a tree that is rooted at the substation. We also define $\mathcal{I}=\{i: i=1,\ldots,I\}$ as the set of poles, $\mathcal{I}^u=\{i: i=1,\ldots,I^u\}$ set of nodes of circuit $u$ which are mounted on the poles  and $\mathfrak{L}^u$ as the set of power lines on circuit $u$ where power line $i\in \mathfrak{L}^u$ feeds node $i\in \mathcal{I}^u$ with power.

Using the configuration of the grid, we define $d^u_i$ to be the first upstream protective device of node $i$ and power lines $i$ on circuit $u$.
For each circuit $u$, let $\mathcal{D}^u$ be the set of protective devices  and $\mathcal{S}^u$ be the set of segments where each segment contains the power lines that trigger the same protective device. We define ${\cal S}^u_i$ as the set of power lines belonging to the same segment of node~$i$ on circuit $u$, i.e., triggering the same protective device, $\mathcal{S}^u_i=\{j, d^u_j = d^u_i \}$. We also introduce $\mathcal{Q}^u_{i}=\{j, \mbox{ node } i\in \mathcal{I}^u \hbox{~becomes in outage if power line~} j \in \mathfrak{L}^u \hbox{~faults}\}$ as the set of power lines that if faulted result in an outage to node $i$ on circuit $u$ and $\mathcal{W}^u_i$ as the set of downstream segments of node~$i$ but with different upstream protective devices, i.e., $\mathcal{W}^u_i=\{S^u_j, \hbox{~segment~}S^u_j \hbox{~is downstream of} \linebreak  \hbox{node~}i ~\hbox{~\&~}  d^u_i\neq d^u_j\}$. We assume that $n_i^u$ customers are attached to node $i$ on circuit $u$. If the node is a transformer, then $n_i^u>0$, otherwise $n^u_i=0$ because no customers are attached to the power generator or protective devices.

Let $G(\mathcal{V},\mathcal{E})$ be the graph representing the road network through which the utility crew can travel to check the power grid where $\mathcal{V}=\{i: i=1,\ldots,N\}$ is the set of nodes of the road network where each node represents either a begin/end of a road segment or a pole of the power grid, and $\mathcal{E}$ is the set containing the arcs/roads in the graph, i.e.,  $\mathcal{E}=\{(i,j): \mbox{~if there is an arc between nodes~} i \mathrm{~and~ } j \mathrm{~in~} \mathcal{V}\}$. Thus, some of the arcs/roads in $\mathcal{E}$ are parallel to the power lines of the power~grid as shown in Figure~\ref{fig:sys_mod}.


In this paper, we assume that an important source of information that the EUCs rely on to identify power outages is the calls of the customers that lost power. Let $H_t=\{H^u_{ti}: \forall i\in \mathcal{I}^u, \forall u\in \mathcal{U}\}$ be a random vector representing the possible realizations of received calls from the nodes of the power system by time $t$ where $H^u_{ti}$ is a random variable representing the number of received phone calls from node $i$ on circuit $u$. Customer phone calls arrive over time, thus, $H_t=H_{t-1}+\hat{H}_t$ where $\hat{H}_t$ represents the new set of received calls in the time interval between $t-1$ and $t$.

Based on the set of received phone calls by time $t$, we have to figure out the power lines that faulted in the grid. For each circuit $u$, let $L^u_t=\{L^u_{ti}: \forall i\in \mathfrak{L}^u\}$ be a random vector representing the possible realizations of power lines that faulted by time $t$  on circuit $u$ where $L^u_{ti}$ is a random variable indicating whether power line $i$ on circuit $u$ has faulted; we assume $L^u_{ti}=1$ if the power line faulted and $L^u_{ti}=0$ otherwise.  

Let $T_{tij}$ be a random variable representing the required travel time from node $i$ to node $j$ which depends on the traffic condition of the road, at time $t$, and $R^u_{i}$ be a random variable representing the required repair time for power line $i$ on circuit $u$ which depends on the fault type. Now, we can define $\tau_{tij}=T_{tij} + \sum_{u\in\mathcal{U}}R^u_{j}$ as the total time required to go from node $i$ to $j$ which accounts for the travel and repair times of all power lines that feed pole $j$ with power.

Let $\omega$ be a sample realization of the random variables. 
At time $t$ and according to sample path $\omega$, let $H_t(\omega)$ be an outcome of customer calls, $L^u_t(\hspace{-0.03cm}\omega \hspace{-0.05cm})$ be an outcome of the power line faults on circuit $u$, $R^u(\omega)$ be an outcome of the repair times on circuit $u$ and $T_{t}(\omega)$ be an outcome of the travel conditions. Then, we can define at time $t$ and according to sample path $\omega$, $p(H_t=H_t(\omega))$ as the probability of the set of received calls, $p(L^u_t=L^u_t(\hspace{-0.03cm}\omega \hspace{-0.05cm}))$ as the prior probability of power line faults on circuit $u$, $p(R^u=R^u(\omega))$ as the probability of the repair time on circuit $u$ and  $p(T_t=T_t(\omega))$ as the probability of the travel time.


Define $\Omega$ as the set of all outcomes, $\cal F$ as the set of events and $\cal P$ as the probability measure on $(\Omega,\cal F)$, so that $(\Omega,\cal F,\cal P)$ is the probability space. $\Omega\subseteq \mathcal{F} $ is formed by a set of scenarios where each scenario~$\omega$ indicates a specific set of calls, a set of power line faults, fault types and travel conditions (traffic). Let $p(\omega)$ be the probability of scenario~$\omega$, where $\sum_{\omega\in\Omega}p(\omega) = 1$.

Second, the EUC should estimate the prior probability, $p(L^u_{ti})$, of fault of power line $i$ on circuit $u$ at time $t$ based on the storm pattern, the structure of the grid and the EUC knowledge of the power line conditions and environment as explained in~(\citet{ABP14}). In the worst case, the EUC can assume a uniform probability of fault for all power lines, and then the set of received of calls will identify the ones that have most likely faulted.

In~(\citet{ABP14}), for any realization of $L^u_t$ and $H_t$, the posterior probability of fault of the power lines on circuit $u$ given the phone calls is calculated using Bayes' theorem as follows:
\begin{eqnarray}
p(L^u_t|H_t) = \frac{p(H_t|L_t^u)p(L_t^u)}{p(H_t)},
\end{eqnarray}
where $p(H_t|L_t^u)$ is the likelihood of the calls given the power line faults on circuit $u$ and the expressions have been derived in equation (5) in~(\citet{ABP14}); it is a function of the locations of the calling customers, the calling probability which refers to the percentage of customers calling to report an outage and the grid structure.  The posterior probability of fault of power line $i$ on circuit $u$, given the phone calls can be expressed as:
\begin{eqnarray}
p(L^u_{ti}=1|H_t)=\frac{\sum_{L^u_t\in \{\mathcal{L}^u\}_{L^u_{ti}=1}}p(H_t,L^u_t)}{p(H_t)}=\frac{\sum_{L^u_t\in \{\mathcal{L}^u\}_{L^u_{ti}=1}}p(H_t|L^u_t)p(L^u_t)}{\sum_{L^u_t\in\mathcal{L}^u} p(H_t|L^u_t)p(L^u_t)}, \label{eq:eq2}
\end{eqnarray}
\noindent where $\mathcal{L}^u$ is the set containing all power line fault combinations on circuit $u$ and  $\{\mathcal{L}^u\}_{L^u_{ti}=1}$ is the set containing a subset of vectors of $\mathcal{L}^u$ where the variable corresponding to power line $i$, i.e., $L^u_{ti}$, is equal to $1$. Thus, $\{\mathcal{L}^u\}_{L^u_{ti}=1}$ is the set containing all the combinations of power lines that can fault with power line $i$ on circuit $u$.

In this work, another factor plays a major role in identifying the fault probability of a power line by time  $t$ which is the trajectory of the truck that is going across the power grid to fix faults. For example, if power line $i$ on circuit $u$ has been fixed by time $t-1$, then its prior probability of fault, at time $t$, is $0$. 
Let $x_{tij}$ be a binary variable representing whether the utility truck travels from node $i$ to node $j$ using roadway graph $G(\mathcal{V},\mathcal{E})$ at time $t$. It is assumed that if a truck travels from node~$i$ to node~$j$  at time $t$ then it repairs all the power lines that are attached to pole $j$ if there are faults across them.


Let $x_t$ be a matrix capturing the vector of decisions $(x_{tij})_{i,j\in\mathcal{V}}$ where $x_{tij} = 1$ if the truck is dispatched to $i$ from $j$ at time $t$. Also, let $X_t$ be the trajectory of the truck up to time $t$, i.e., $X_t = (x_{t'})_{t'=0}^t$. The information we are looking for, at time $t$, is the posterior probability, $p(L^u_{ti}|H_t,X_{t-1})$,  of power line $i$ being in fault given the phone calls and the trajectory of the truck up to time $t-1$ which is given by:
\begin{eqnarray}
&&\hspace{-0.7cm}p(L^u_{ti}=1|H_t,X_{t-1})=\frac{\sum_{L^u_t\in \{\mathcal{L}^u\}_{L^u_{ti}=1}}p(H_t|L^u_t)p(L^u_t|X_{t-1})}{\sum_{L^u_t\in\mathcal{L}^u}p(H_t|L^u_t)p(L^u_t|X_{t-1})}\label{eq:eq3},
\end{eqnarray}
where $p(L^u_t|X_{t-1})$ is the prior probability of vector $L^u_t$ being in fault given the route of the truck; the prior probability of the power lines is updated by setting $p(L^u_{ti}|X_{t-1})=0$ if $\sum_j x_{t'ji}=1, t'\leq t-1$. However, the likelihood $p(H_t|L^u_t)$ is independent of the route of the truck. Accordingly, at each time $t$, the prior probabilities of power line faults should be computed and the set of received calls should account for any new incoming calls between $t-1$ and $t$. After collecting the updated priors and customer calls, the posterior probability of faults are computed using the model described in~(\citet{ABP14}) upon which the route of the truck is determined, i.e., the value of $x_{t}$.

The aim of this work is to develop a policy that routes a utility truck to restore outages in the power grid as quickly as possible while accounting for all the components described above.

\section{Sequential Stochastic Optimization Problem}~\label{sec:opt_prob}
To develop a near-optimal policy for managing a utility truck across the grid, the problem should be formulated first as a sequential stochastic optimization problem after which we can derive the optimal policy to solve it.  However, this introduces several challenges as explained below.

First, the number of customers whose power is restored after a truck visits a location (even if a repair occurs) is a random variable since the actual value depends also on upstream and downstream outages which are uncertain. Second, each time a truck travels a segment, the information collected (e.g. that there is an outage on that segment, or not) is used to update the probability of outages on all lines. Third, new phone calls are arriving over time, which allows us to update probabilities of outages.  At the same time, as trucks identify and fix outages, other phone calls may become irrelevant. Finally, the time required for a truck to traverse a segment depends to a large extent on whether it finds an outage, and the time required to repair the outage.

\begin{figure}[t!]
\vspace{-0cm}
\centering
\includegraphics[scale=0.6]{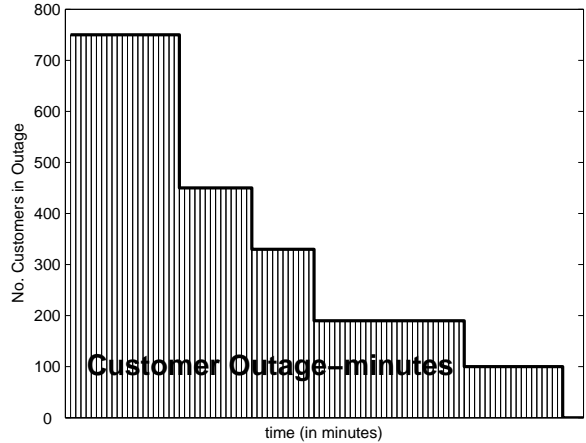}
\vspace{-0.3cm}
\caption{\label{fig:obj}Objective function; outage-minute is represented by the shaded area under the curve.}\vspace{-0.3cm}
\end{figure}

In order to develop the sequential stochastic optimization problem that restores power to the maximum number of customers in outage over time, we need to define its five fundamental elements:
\begin{itemize}
\item State $S_t$ - The information capturing what we know at time $t$; in this work, $S_t=\bigg(R_t,P_t^L,H_t\bigg)$ where $R_t$ represents the physical state that indicates the location of the truck at time $t$, $P_t^L$ is a vector containing the prior probability of all power line faults (i.e., its entries are $p(L_{ti}^u=1|X_{t-1})$ and $H_t$ is the set of received calls by time $t$). Thus, given the prior probabilities of fault and the set of calls, we can calculate the posterior probabilities of faults which also represent the state of the network at time $t$.
\item Decision $x_t$ - The vector $x_t = (x_{tij})_{i,j}$  captures the decision made at time $t$ about the next hop of the truck, where $x_{tij} = 1$ if the truck is sent from node $i$ to node $j$ at time $t$. Let $X^\pi(S_t)$ be the policy that determines $x_t \in \mathcal{X}_t$ given~$S_t$.
\item Exogenous information $W_t$ - The new information that arrives between $t-1$ and $t$.  This includes new phone calls ${\hat H}_t$, as well as information about outages (discovered by the utility trucks) and travel (or repair) times.  We denote the exogenous information process by $W_1, W_2, ..., W_T$.  Let $\omega$ be a sample realization of the information, which we note depends on the policy (if a truck fixes an outage, then this will produce fewer phone calls).  For this reason, we let $\Omega^\pi$ be the set of outcomes which depends on the policy $\pi$ that we use to dispatch trucks. Let ${\cal F}^\pi$ be the sigma-algebra on $\Omega^\pi$, and let ${\cal P}^\pi$ be the probability measure on $(\Omega^\pi, {\cal F}^\pi)$, giving us a probability space $(\Omega^\pi, {\cal F}^\pi, {\cal P}^\pi)$.  We let ${\cal F}^\pi_t$ be the sigma-algebra generated by $W_1, ..., W_t$, giving us the filtrations ${\cal F}^\pi_t \subseteq {\cal F}^\pi_{t+1}$.  This notation means that all variables indexed by $t$ are ${\cal F}^\pi_t$-measurable.

\item The transition function $S_{t+1}=S^M(S_t,x_t,W_{t+1})$ which represents the evolution of the physical, informational and belief states; for example, if power line $j$ on circuit $u$ was repaired at time $t$ then its prior probability of fault is set to $0$ at time $t+1$, i.e., $p(L_{{t+1}j}=1|\sum_ix_{tij}=1)=0$. This in turn updates  the probabilities of faults over the power grid. The function also represents the exogenous evolution of calls, i.e., $H_{t+1}=H_t+\hat{H}_{t+1}$ and any additional information contained by $W_t$ as explained earlier.
\item The objective function - Our objective is to find the policy that maximizes the number of customers with restored power over time which is equivalent to minimizing the number of customers in outage while minimizing the truck's operating cost. Let $C(S_t,x_t) = C^c(S_t,x_t)+\gamma C^o(S_t,x_t)$ be the total cost of taking decision $x_t$ given state $S_t$ where $C^c(S_t,x_t)$ and $C^o(S_t,x_t)$ represent the customer outage-minutes and the truck's operating costs, respectively, and $\gamma$ is a tunable parameter that balances the weight between the two costs.  The cost function $C^c(S_t,x_t)$ represents the number of customers in outage at time $t$ as shown in Figure~\ref{fig:obj} and the sum over time, $\sum_{t=0}^TC^c(S_t,x_t)$, evaluates the shaded area under the curve which we refer to as ``customer outage-minutes."
    The total objective can be represented as $\sum_{t=0}^TC(S_t,x_t)$ and we have to find the policy that solves:
    \begin{eqnarray}
    \min_\pi \mathbb{E}^\pi\left[\sum_{t=0}^TC(S_t,X^{\pi}(S_t))|S_0\right]\label{eq:optimal_policy}
    \end{eqnarray}
    where the expectation is over all possible sequences $W_1,W_2,\ldots,W_T$, which  depend on the decisions taken. The initial state $S_0$ captures all deterministic parameters, and priors on any uncertain information (such as outages). Our goal is to find the best policy for making a decision.
\end{itemize}

So far, we have defined the fundamental elements influencing truck routing. In this work, the flow of information and state transition occurs in the following sequence:
\begin{eqnarray}
(S_0,x_0,S^x_0,W_1,S_1,x_1,S^x_1,W_2,S_2,\ldots,S_t,x_t,S^x_t,W_t,\ldots,S_T)
\end{eqnarray}
where $S^x_t$ represents the state after taking decision $x_t$, known as the post-decision state. Based on $x_t$, we see a specific exogenous information $W_{t+1}$ indicating that it directly depends on~$x_t$.

The distinguishing feature of this problem relative to the classical stochastic vehicle routing literature is the dimension that the utility truck is collecting information, and that our state variable includes the physical state of the truck, our current state of knowledge concerning the probability of outages, and the history of phone calls (there is a reason we have to retain the history).  A decision to dispatch a truck from $i$ to $j$ has to consider not only the change in the physical state of the truck but also the value of the information that is collected while traversing from $i$ to $j$.

The contribution of this paper relative to the dynamic vehicle routing literature is that we capture the collection of information (phone calls, observations of outages) to update our beliefs about where outages may be.  In contrast with classical models where the state $S_t$ captures the status of trucks and unserved demands, our model uses a state variable with three types of information:  the physical state $R_t$, which includes the location of the truck, information $I_t$ which includes the history of phone calls $H_t$ and observations of actual outages, and the probabilistic knowledge $P^L_t$ which captures the probability of outages on observed lines.  Thus we write the state of our system as $S_t = (R_t, H_t, P^L_t)$.  We update our probabilistic knowledge $P^L_t$ by applying Bayes' theorem that combines the prior $P^L_t$ with exogenous information of phone calls and grid observations to produce the posterior~$P^L_{t+1}$.


\citet{Powell14} describes four classes of policies for solving the optimization in equation (\ref{eq:optimal_policy}): policy function approximations (PFAs), parametric cost function approximations (CFAs), policies based on value function approximations (VFAs), and lookahead policies.  The utility currently uses a rule-based PFA (send a truck to the next caller reporting an outage).  In this paper, we are going to use a policy based on a stochastic lookahead~model.

\subsection{Multistage Lookahead Policy}

A lookahead policy is particularly useful for our problem since it is time-dependent obviously because the grid has a radial structure which affects the computation of the objective.

First, we begin our discussion by characterizing an optimal policy using the function
\begin{eqnarray}
X^*_t(S_t)=\arg\min_{x_t\in{\mathcal{X}_t(S_t)}}\left(C(S_t,x_t)+\min_\pi \mathbb{E}^\pi\left\{ \sum_{t'=t+1}^T C(S_{t'},X^\pi_{t'}(S_{t'}))|S^x_t\right\} \right), \label{eq:opt_policy}
\end{eqnarray}
where $S_{t+1}=S^M\left(S_t, x_t,W_{t+1}\right)$. Generally, (\ref{eq:opt_policy}) represents a lookahead policy solving a multistage stochastic program , i.e., after implementing the first-period decision and stepping forward in time, the problem has to be solved again with a new set of scenario trees. In scenario trees, we choose $x_t$, then see the information $W_{t+1}$, then we would choose $x_{t+1}$, after which we see $W_{t+2}$, and so on. Thus, (\ref{eq:opt_policy}) is equivalent to:
\begin{eqnarray}
X^*_t(S_t)&=&\arg\min_{x_t\in{\mathcal{X}_t(S_t)}}\bigg(C(S_t,x_t)+\mathbb{E}_{W_{t+1}\in\Omega_{t+1}}\bigg[\min_{x_{t+1}\in{\mathcal{X}_{t+1}(S_{t+1})}}C(S_{t+1},x_{t+1})+\mathbb{E}_{W_{t+2}\in\Omega_{t+2}}\bigg[ \ldots + \nonumber\\
&&\hspace{5.5cm}\mathbb{E}_{W_{T}\in\Omega_{T}}\bigg[C(S_{T})|S^x_{T-1}\bigg]\ldots\bigg]|S^x_{t+1}\bigg]|S^x_t \bigg]\bigg).\label{eq:eqv_opt_policy}
\end{eqnarray}
Scenario trees capture the entire history; this is particularly important for our problem since the history of the route of the truck which indicates the power lines that have been fixed affects the computation of the customers in outage due to the radial structure of the grid.

In practice, the policies in (\ref{eq:opt_policy}) or (\ref{eq:eqv_opt_policy}) are computationally intractable, requiring that we introduce an approximation. There are several strategies that are typically used to simplify lookahead models:
 \begin{itemize}
\item[1) ] Limiting the horizon by reducing it from $(t,T)$ to $(t,t+H)$.
\item[2) ] Discretizing the time, states and decisions to make the model computationally tractable.
\item [3) ] Aggregating the outcome or sampling by using Monte Carlo sampling to choose a small set, $\tilde{\Omega}_t$, of possible outcomes between $t$ and $t+H$.
\item[4)] Stage aggregation which represents the process of revealing information before making another decision; a common approximation is a two-stage formulation, where we make
a decision $x_t$, then observe all future events (until $t+H$), and then make all remaining decisions. A more accurate formulation is a multistage model but these can be computationally expensive if not strategically designed.
 \item[5) ] Dimensionality reduction where we ignore some variables in our lookahead model as a form of simplification. For example, a forecast of future incoming phone calls can add a number of dimensions to the state variable. While we have to track these in the original model, we can hold them fixed in the lookahead model, and then ignore them in the state variable (these become latent variables).
\end{itemize}

The most common lookahead strategy is a deterministic model, although even these can be difficult for vehicle routing problems.  For this reason~(\citet{GoOh2013}) explores the use of a rollout policy based on ideas from~(\citet{BeCa}) which uses a simple heuristic to approximate the downstream value of a decision.  Deterministic lookahead models are unable to capture the uncertainty of our belief state, while rollout heuristics, aside from offering no theoretical guarantees, are typically built from myopic policies which are unlikely to work well.

All variables in the lookahead model are indexed by $(t,t')$ where $t$ represents when the decision is being made (which fixes the information content) while $t'$ is the time within the lookahead model. We also use tilde's to avoid confusion between the lookahead model (which often uses a variety of approximations) and the real model. So, the stochastic lookahead model becomes
\begin{eqnarray}
\hspace{-0.5cm}X^*_t(S_t)&=&\arg\min_{x_t\in{\mathcal{X}_t(S_t)}}\bigg(C(S_t,x_t)+\tilde{\mathbb{E}}_{\tilde{W}_{t,t+1}\in\tilde{\Omega}_{t,t+1}}\bigg[\min_{\tilde{x}_{t,t+1}\in{\mathcal{\tilde{X}}_{t,t+1}(\tilde{S}_{t,t+1})}}\tilde{C}(\tilde{S}_{t,t+1},\tilde{x}_{t,t+1})+ \nonumber\\
&&\hspace{-0.7cm}\mathbb{\tilde{E}}_{\tilde{W}_{t,t+2}\in\tilde{\Omega}_{t,t+2}}\bigg[ \ldots\mathbb{\tilde{E}}_{\tilde{W}_{t,t+H}\in\tilde{\Omega}_{t,t+H}}\bigg[\tilde{C}(\tilde{S}_{t,t+H})|\tilde{S}^x_{t,t+H-1}\bigg]\ldots\bigg]|\tilde{S}^x_{t,t+1}\bigg]|S^x_t \bigg]\bigg),\label{eq:alo}
\end{eqnarray}
where the expectation $\tilde{\mathbb{E}}\{.|S^x_t\}$ is over the sample space in $\tilde{\Omega}_{t,t+1}$ which is constructed given that
we are in state $S_t$ at time $t$. When computing this policy, we start in a particular state $S_t$, but then step forward in time using:
\begin{eqnarray}
\tilde{S}_{t,t'+1}=S^M(\tilde{S}_{tt'},\tilde{x}_{tt'},\tilde{W}_{t,t'+1}), t'=t,\ldots,t+H-1.
\end{eqnarray}
Alternatively, we can formulate~(\ref{eq:alo}) as:
\begin{eqnarray}
\hspace{-0.5cm}X^*_t(S_t)&=&\arg\min_{x_t\in{\mathcal{X}_t(S_t)}}\bigg(C(S_t,x_t)+\tilde{V}^x_t(S^x_t)\bigg),\label{eq:al}
\end{eqnarray}
where $\tilde{V}^x_t(S^x_t)$ is the approximate average value of the \emph{post-decision state} that is returned by the lookahead policy.

In this work, we fix the set of received calls at time $t$ while the other sources of randomness such as the fault locations and the travel/repair times are included in the tree. This lookahead model could be solved via scenario trees but the size of the tree can grow exponentially if we take into account all possible random realizations since taking decision $\tilde{x}_{tt'}$ from state $\tilde{S}_{tt'}$ can result in many different successor states based on the observed exogenous information. To simplify the problem, we can assume that the only randomness that is modeled in the tree is whether a fault is found at a certain location or not. Also, one fault type can be considered which results in a deterministic repair time and the travel times can be also considered deterministic but they definitely depend on the distances between the nodes in the grid. In this case, there are only two possible successor states by taking decision $x_{tt'}$ from state $S_{tt'}$ which depends on whether a fault is observed or not.


The number of customers whose power is restored due to decision $x_t$ from state $S_t$ is a random variable given the uncertainty about the location of outages due to the radial structure of the grid. The actual number of served customers depends on whether there is a fault upstream to the visited location, which is uncertain. In addition, if a location faults, it causes outage to the customers attached to its segment and all downstream segments. Thus, if a fault is fixed, then all the downstream customers might restore power unless there are other downstream segments that faulted which is uncertain. The actual ``customer outage-minutes" objective is represented in Figure~\ref{fig:obj}, but since we cannot figure out the exact number of customers with restored power at each state $S_t$, we evaluate the expected customer outage-minutes. At state $S_t$, the fault probabilities are updated according to (\ref{eq:eq3}) and the expected number of customers in outage~is:
    \begin{eqnarray}
    C^c(S_t,x_t)=   \sum_{u\in\mathcal{U}}\sum_{s\in S^u} \left(1-\prod_{k\in Q_s}p(L^u_{tk}=0)\right)\sum_{k\in s} n^u_k,
 \label{eq:ecom}
    \end{eqnarray}
where $Q_s$ represents the set of power lines that if fault result in outage to segment $s$, the term in parenthesis is the probability of at least one fault across $Q_s$ that causes $s$ to be in outage and $\sum_{k\in s} n^u_k$ is the number of customers across segment $s$.

Generating the entire scenario tree is computationally intractable.  For this reason, we turn to a popular strategy developed in the computer science community known as Monte Carlo tree search (MCTS)~(\citet{Munos14,BPW12,CHS08,KS06}) which uses an intelligent sampling procedure to create a partial tree that can be solved, and which asymptotically ensures that we would find the optimal decision given a large-enough computing budget.  Given that we are solving the one-truck problem, we can formulate the lookahead model as a decision tree, taking advantage of the property that the number of possible decisions for a truck at any point in time is quite small, since it is limited to paths over a road network.  Further, if we limit the random information to whether the truck finds a fault or not, then the random variables are binomial.  However, even with these restrictions, a decision tree will still grow exponentially. 


We propose  MCTS as the look-ahead policy to solve the original problem in the following way. Given a current state, $S_t$, which depends on the location of the truck and the probability of faults at time~$t$,  MCTS should decide where to move the truck next in order to minimize the objective represented by the customer outage-minutes. So, starting from the current state $S_t$ as the root node, sampled MCTS successively builds the look-ahead tree over the state-space. Finally, the move that corresponds to the highest value from the root node will be taken. At this stage, the probability space is updated taking into account any incoming exogenous information such as whether there was a fault or not by taking the move that was the outcome of the previous step, the newly arrived phone calls of customers and the consumed travel/repair times. Then, the whole process is repeated until power is reconstructed to the whole power grid based on the values of the probability model; that is, when the values of the posterior probabilities drop below a certain threshold $\epsilon^{thr}$ which is typically a small value. The pseudo-code of the proposed lookahead policy to solve the utility truck routing problem is presented in Algorithm~\ref{alg:adp}.


\begin{algorithm}[t!]{\footnotesize
\caption{Lookahead Policy for Utility Truck Routing}
\label{alg:adp}
\begin{algorithmic}[t!]
\STATE\textbf{Step 0.} \textbf{Initialization}:  Initialize the state $S_0=(R_0,P_0^L,H_0)$. Set $t\leftarrow0$.
\STATE\textbf{Step 1.}~\textbf{While} $\left(p(L^u_{t{i}}|H_t,X_{t-1})\geq \epsilon^{thr}, \forall i,\forall u\right)$ do:
\STATE \hspace{1.2cm}\textbf{Step 1a.} At time $t$, fix the set of calls $H_t$ and determine $\tilde{\Omega}_t$.
\STATE \hspace{1.2cm}\textbf{Step 1b.} Call $MCTS(S_t)$ to solve (\ref{eq:al}) that determines $x^*_t$.
\STATE \hspace{1.2cm}\textbf{Step 1c.} Set $x_t\leftarrow x^*_t$ and move the truck according to $x_t$.
\STATE \hspace{1.2cm}\textbf{Step 1d.} Update $S_{t+1}=S^M(S_t,x_t,W_{t+1})$ until the truck reaches the destination node set by $x^*_t$, say at \STATE \hspace{2.5cm} time $t'$.
\STATE \hspace{1.2cm}\textbf{Step 1e.} $t\leftarrow t'$
\STATE \hspace{1cm}\textbf{End While}
\end{algorithmic}\vspace{-0.0cm}
}
\end{algorithm}

\section{Monte Carlo Tree Search}\label{sec:MCTS}

Even if we are solving a stochastic lookahead for a single truck, the resulting decision tree still grows exponentially, and is too large to enumerate.  For this reason, we turn to a technique known as Monte Carlo tree search (MCTS), widely used in computer science, which dynamically creates a decision tree using a combination of Monte Carlo sampling and heuristic rollout policies to identify the most promising branches.

MCTS builds successively a search tree until some predefined computational budget is reached such as number of iterations or time constraints. In deterministic MCTS, each node corresponds only to the physical state of the network and each edge corresponds to a decision which results in another unique physical state. In stochastic MCTS,  the exogenous information is included in the tree and thus each node corresponds to a state and each edge can represent either a decision (move) or an exogenous event both of which result in another state. Each iteration of MCTS is formed of four basic steps which are Selection, Expansion, Simulation and Backpropagation~(\citet{CHS08}). The \emph{Selection} step involves selecting a decision successively starting from the initial state till an expandable state is reached.  The \emph{Expansion} step adds one or more states to the tree. The \emph{Simulation} step referred to as ``Simulation Policy" is used to evaluate the value of the newly added state. Finally, in the \emph{Backpropagation} step, the value of the newly added state is backpropagated to update the value functions of all predecessor states.

MCTS has evolved in the literature primarily in the context of deterministic problems. Deterministic MCTS has been used extensively in the literature, as summarized in~(\citet{BPW12}). In this case, the exogenous information is handled in the simulation policy by averaging over the possible outcomes. Stochastic outcomes (which characterizes our application) has been handled using a process known in the computer science community as ``{\it Determinization}"~(\citet{BFT09,BSC07,Caz06}), or by explicitly representing exogenous information in the tree directly~(\citet{CHS11}), which is the approach that we take. It is shown in~(\citet{ACT13}) that stochastic MCTS, in which the exogenous information is included in the tree, is asymptotically optimal if every possible action can be uniformly chosen to be included in the tree no matter how good or bad it is. The same rule applies for the exogenous events. In this work, we choose a more strategic way in choosing the actions to expand in order to converge faster to the optimal solution.



In stochastic MCTS, which is also referred to as sampled MCTS, the states in the tree are associated with the following~data:
 \begin{itemize}
 \item[1) ] The pre-decision value function, $\tilde{V}_{tt'}(\tilde{S}_{tt'})$, the post-decision value function, $\tilde{V}^x_{tt'}(\tilde{S}^x_{tt'})$, and cost, $\tilde{C}(\tilde{S}_{tt'},\tilde{x}_{tt'})$, which represent the value function and cost of being in state $\tilde{S}_{tt'}$ and taking decision $\tilde{x}_{tt'}$, respectively.
 \item[2) ] The visit count, $N(\tilde{S}_{tt'})$, which represents the number of rollouts that included state $\tilde{S}_{tt'}$.
 \item[3) ] The count of the state-decision, $N(\tilde{S}_{tt'},\tilde{x}_{tt'})$, which represents the number of times decision $\tilde{x}_{tt'}$ was taken from state $\tilde{S}_{tt'}$.
 \item[4) ] The set of decisions, $\tilde{\mathcal{X}}_{tt'}(\tilde{S}_{tt'})$, and set of possible random outcomes, $\tilde{\Omega}_{t,t'+1}(\tilde{S}_{tt'}^{x})$; $\tilde{\mathcal{X}}_{tt'}(\tilde{S}_{tt'})$ is the set of decisions that the truck would face moving over a road network given that it is at state $\tilde{S}_{tt'}$. For state $\tilde{S}_{tt'}$, let $\tilde{\mathcal{X}}_{tt'}^e(\tilde{S}_{tt'})$ be the set of decisions that has been explored by the truck (i.e., expanded in the tree)  by time $t'$ and let $\tilde{\mathcal{X}}_{tt'}^u(\tilde{S}_{tt'})$ be its complement set which represents the set of unexplored decisions in the tree by time $t'$. Similarly,  $\tilde{\Omega}_{t,t'+1}(\tilde{S}_{tt'}^{x})$ is the set of all possible random events that can take place at time $t'+1$ given state $\tilde{S}_{tt'}^{x}$, $\tilde{\Omega}_{t,t'+1}^e(\tilde{S}_{tt'}^{x})$ is the set of explored events and $\tilde{\Omega}_{t,t'+1}^u(\tilde{S}_{tt'}^{x})$ is its complement.
 \end{itemize}

In sampled MCTS, assume that we are at a \emph{pre-decision state} $\tilde{S}_{tt'}$ and decide to take decision $\tilde{x}_{tt'}$ which takes us to the \emph{post-decision state} $\tilde{S}_{tt'}^{x}$. In real life, while the truck is moving to the location specified by $x_{t}$, it encounters first the travel time which depends on traffic. Second,  it discovers the fault type at the intended location which determines the repair time and affects the number of customers with restored power. In the lookahead model, upon determining $\tilde{x}_{tt'}$, we sample one of the possible exogenous realizations $\tilde{W}_{t,t'+1}$ which immediately informs us about the expected travel and repair times and thus we can immediately know the time, $\tau(\tilde{x}_{tt'},\tilde{W}_{t,t'+1})$, required by the truck to arrive to the destination node specified by $\tilde{x}_{tt'}$. Thus, in MCTS, we define the stochastic transition function as $\tilde{S}_{t,t'+\tau(\tilde{x}_{tt'},\tilde{W}_{t,t'+1})}=\tilde{S}^{M,x}(\tilde{S}_{tt'}^{x},\tilde{W}_{t,t'+1})$. It is obvious now that taking the same decision $\tilde{x}_{tt'}$ from the same state $\tilde{S}_{tt'}$ results in a different outcome state based on the exogenous information $\tilde{W}_{t,t'+1}$.

\begin{figure}[t!]
\vspace{-0cm}
\centering
\includegraphics[scale=0.36]{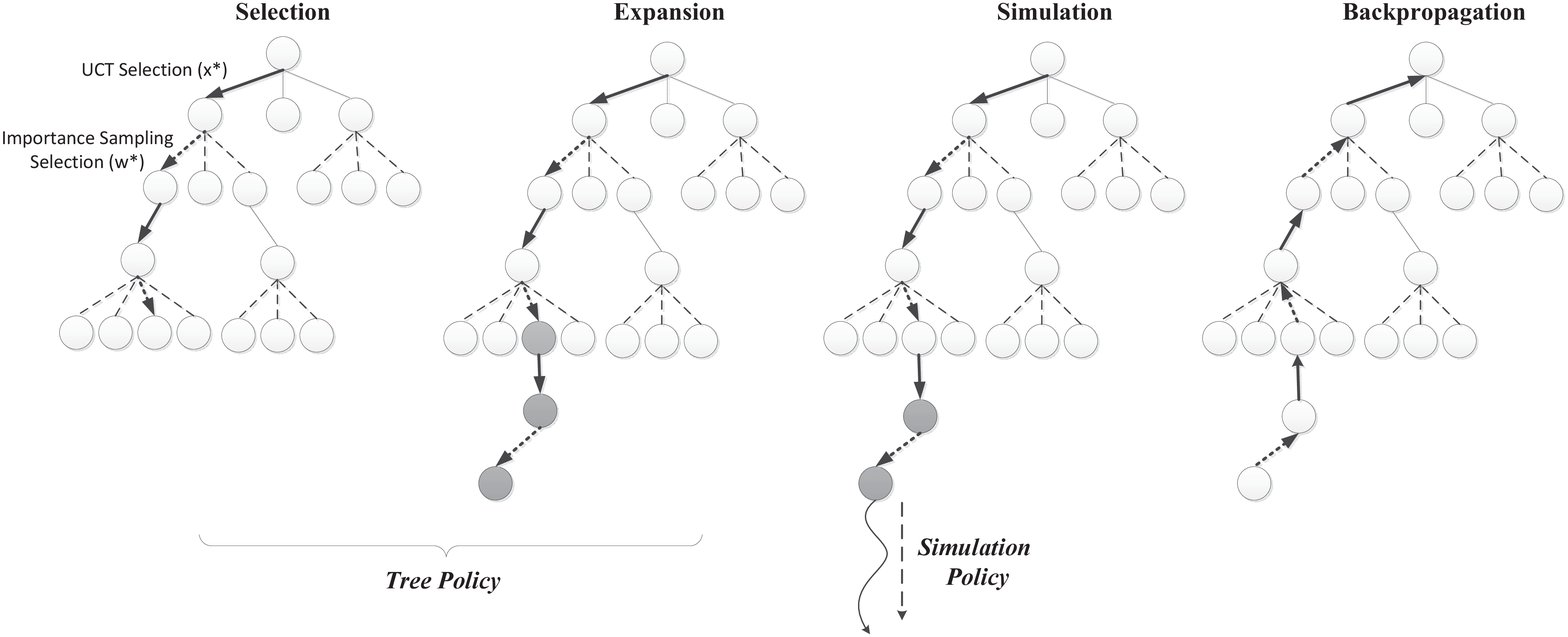}
\vspace{-0.3cm}
\caption{\label{fig:S_MCTS}One iteration of the proposed sampled MCTS.}\vspace{-0.3cm}
\end{figure}

We assume that MCTS has a computational budget of $n^{thr}$ iterations. The steps of MCTS can be grouped into two main policies: a \emph{Tree Policy } (formed of the \emph{Selection} and \emph{Expansion} steps) and a \emph{Simulation Policy} (formed of the \emph{Simulation} step) as shown in Figure~\ref{fig:S_MCTS}. After terminating the tree search, then the decision that corresponds to best value from the root node is chosen. The pseudo-code for sampled MCTS is presented in Algorithm~\ref{alg:MCTS}.



\begin{algorithm}[t!] \footnotesize
\caption{Sampled MCTS Algorithm}
\label{alg:MCTS}
\begin{algorithmic}[t!]
\STATE \textbf{function} $MCTS(S_t)$
\STATE \hspace{0.5cm}Create root node $\tilde{S}_{tt}$ with state $S_t$; set iteration counter $n=0$
\STATE \hspace{0.5cm}{\bf while} $n< n^{thr}$
\STATE \hspace{1cm} $\tilde{S}_{tt'}\leftarrow TreePolicy(\tilde{S}_{tt})$
\STATE \hspace{1cm} $\tilde{V}_{tt'}(\tilde{S}_{tt'})\leftarrow SimPolicy(\tilde{S}_{tt'})$
\STATE \hspace{1cm} $Backup(\tilde{S}_{tt'},\tilde{V}_{tt'}(\tilde{S}_{tt'}))$
\STATE \hspace{1cm} $n\leftarrow n+1$
\STATE \hspace{0.5cm}{\bf end while}
\STATE \textbf{return} $x_{t}^*=\arg\min_{\tilde{x}_{tt}\in \tilde{\mathcal{X}}_{tt}^e(\tilde{S}_{tt})}\tilde{C}(\tilde{S}_{tt},\tilde{x}_{tt})+\tilde{V}^x_{tt}(\tilde{S}^x_{tt})$
\STATE
\STATE \textbf{function} $TreePolicy(\tilde{S}_{tt})$
\STATE \hspace{0.5cm} $t'\leftarrow t$
\STATE \textbf{while} $\tilde{S}_{tt'}$  is non-terminal \textbf{do}
\STATE \hspace{0.2cm} \textbf{if} $|\tilde{\mathcal{X}}_{tt'}^e(\tilde{S}_{tt'})|<d^{thr}$ \textbf{do}(Expanding a decision out of a pre-decision state)
\STATE \hspace{0.5cm} choose decision  ${\tilde{x}_{tt'}}^*$ optimistically by using a two-stage lookahead model
\STATE \hspace{0.5cm} $\tilde{S}_{tt'}^{x}=S^M(\tilde{S}_{tt'},{\tilde{x}_{tt'}}^*)$ (Expansion step)
\STATE \hspace{0.5cm}  $\tilde{\mathcal{X}}_{tt'}^e(\tilde{S}_{tt'})\leftarrow \tilde{\mathcal{X}}_{tt'}^e(\tilde{S}_{tt'})\bigcup\{\tilde{x}_{tt'}^*\}$
\STATE \hspace{0.5cm}  $\tilde{\mathcal{X}}_{tt'}^u(\tilde{S}_{tt'})\leftarrow \tilde{\mathcal{X}}_{tt'}^u(\tilde{S}_{tt'})-\{\tilde{x}_{tt'}^*\}$
\STATE \hspace{0.2cm} \textbf{else}
\STATE \hspace{0.5cm} $\tilde{x}_{tt'}^*=\arg\max_{\tilde{x}_{tt'}\in \tilde{\mathcal{X}}_{tt'}^e(\tilde{S}_{tt'})}\left(-\left(\tilde{C}(\tilde{S}_{tt'},\tilde{x}_{tt'})+\tilde{V}^x_{tt'}(\tilde{S}^x_{tt'})\right)+\alpha\sqrt{\frac{\ln N(\tilde{S}_{tt'})}{N(\tilde{S}_{tt'},\tilde{x}_{tt'})}}\right)$
\STATE \hspace{0.5cm} $\tilde{S}_{tt'}^{x}=S^M(\tilde{S}_{tt'},{\tilde{x}_{tt'}}^*)$
\STATE \hspace{0.2cm} \textbf{end if}
\STATE \hspace{0.2cm} \textbf{if} $|\tilde{\Omega}_{t,t'+1}^e(\tilde{S}_{tt'}^{x})|<e^{thr}$  \textbf{do} (Expanding an exogenous outcome out of a post-decision state)
\STATE \hspace{0.5cm} choose exogenous event ${\tilde{W}_{t,t'+1}}$ according to importance sampling with uniform distribution $g(\tilde{W}_{t,t'+1})$,\STATE \hspace{13cm} $\forall \tilde{\omega}\in\tilde{\Omega}_{t,t'+1}^u(\tilde{S}_{tt'}^{x})$
\STATE \hspace{0.5cm} $\tilde{S}_{t,t'+\tau(\tilde{x}_{tt'},\tilde{W}_{t,t'+1})}=S^{M,x}(\tilde{S}_{tt'}^{x},\tilde{W}_{t,t'+1})$ (Expansion step)
\STATE \hspace{0.5cm} $\tilde{\Omega}_{t,t'+1}^e(\tilde{S}_{tt'}^{x}) \leftarrow\tilde{\Omega}_{t,t'+1}^e(\tilde{S}_{tt'}^{x})\bigcup\{\tilde{W}_{t,t'+1}\}$
\STATE \hspace{0.5cm} $\tilde{\Omega}_{t,t'+1}^u(\tilde{S}_{tt'}^{x}) \leftarrow\tilde{\Omega}_{t,t'+1}^u(\tilde{S}_{tt'}^{x})-\{\tilde{W}_{t,t'+1}\}$
\STATE \hspace{0.5cm} $t'\leftarrow t'+\tau(\tilde{x}_{tt'},\tilde{W}_{t,t'+1})$
\STATE \hspace{0.5cm} \textbf{return} $\tilde{S}_{tt'}$ (stops execution of \textbf{while} loop)
\STATE \hspace{0.2cm} \textbf{else}
\STATE \hspace{0.5cm} choose exogenous event ${\tilde{W}_{t,t'+1}}$ according to importance sampling with uniform distribution $g(\tilde{W}_{t,t'+1})$,\STATE \hspace{13cm}  $\forall \tilde{\omega}\in\tilde{\Omega}_{t,t'+1}^e(\tilde{S}_{tt'}^{x})$
\STATE \hspace{0.5cm} $\tilde{S}_{t,t'+\tau(\tilde{x}_{tt'},\tilde{W}_{t,t'+1})}=S^{M,x}(\tilde{S}_{tt'}^{x},\tilde{W}_{t,t'+1})$
\STATE \hspace{0.5cm} $t'\leftarrow t'+\tau(\tilde{x}_{tt'},\tilde{W}_{t,t'+1})$
\STATE \hspace{0.2cm} \textbf{end if}
\STATE \textbf{end while}
\STATE
\STATE \textbf{function} $SimPolicy(\tilde{S}_{tt'})$
\STATE \hspace{0.5cm} Choose a sample path $\tilde{\omega}\in \tilde{\Omega}_{tt'}$
\STATE \hspace{0.5cm}\textbf{while} $\tilde{S}_{tt'}$ is non-terminal
\STATE \hspace{1cm}Choose $\tilde{x}_{tt'}\leftarrow \pi_0({\tilde{S}_{tt'}})$
\end{algorithmic}
\end{algorithm}

\begin{algorithm}[t!] \footnotesize
\begin{algorithmic}[t!]
\STATE \hspace{1cm} ${\tilde{S}_{t,t'+\tau(\tilde{x}_{tt'}(\tilde{\omega}))}}\leftarrow S^M(\tilde{S}_{tt'},\tilde{x}_{tt'}(\tilde{\omega}))$
\STATE \hspace{1cm} $t'\leftarrow t'+\tau(\tilde{x}_{tt'}(\tilde{\omega}))$
\STATE \hspace{0.5cm}\textbf{end while}
\STATE \textbf{return }$\tilde{V}_{tt'}(\tilde{S}_{tt'})$ (Value function of $\tilde{S}_{tt'}$)
\STATE \textbf{function} $Backup(\tilde{S}_{tt'},\tilde{V}_{tt'}(\tilde{S}_{tt'}))$
\STATE \hspace{0.5cm}\textbf{while} $\tilde{S}_{tt'}$ is not null\textbf{ do}
\STATE \hspace{1cm} $N(\tilde{S}_{tt'})\leftarrow N(\tilde{S}_{tt'}) + 1$
\STATE \hspace{1cm} $t^*\leftarrow$ time when the truck was at predecessor node, i.e., $(\tilde{S}_{tt'}=S^{M,x}(\tilde{S}_{tt^*}^x,{\tilde{W}_{t,t^*+1}}))$ where \STATE \hspace{12cm}$(\tilde{S}_{tt^*}^x=S^{M}(\tilde{S}_{tt^*},{\tilde{x}_{tt^*}}))$
\STATE \hspace{1cm} $\tilde{S}_{tt^*}^{x}\leftarrow \mathrm{~predecessor~of~} \tilde{S}_{tt'}$
\STATE \hspace{1cm} $N(\tilde{S}_{tt^*},{\tilde{x}_{tt^*}})\leftarrow N(\tilde{S}_{tt^*},{\tilde{x}_{tt^*}}) + 1$
\STATE \hspace{1cm} $\tilde{V}^x_{tt^*}(\tilde{S}^x_{tt^*})\leftarrow \frac{1}{\sum_{\tilde{W}_{t,t^*+1}\in\tilde{\Omega}_{t,t^*+1}^e(\tilde{S}_{tt^*}^{x})} p(\tilde{W}_{t,t^*+1})} \cdot E_{g}[p(\tilde{W}_{t,t^*+1})/g(\tilde{W}_{t,t^*+1})\tilde{V}_{tt^*}(S^{M,x}(\tilde{S}_{tt^*}^x,{\tilde{W}_{t,t^*+1}}))]$
\STATE \hspace{1cm} $\tilde{S}_{tt^*}\leftarrow \mathrm{~predecessor~of~} \tilde{S}_{tt^*}^{x}$
\STATE \hspace{1cm} $\Delta \leftarrow \tilde{C}(\tilde{S}_{tt^*},{\tilde{x}_{tt^*}}) + \tilde{V}^x_{tt^*}(\tilde{S}^x_{tt^*}) $
\STATE  \hspace{1cm}  $\tilde{V}_{tt^*}(\tilde{S}_{tt^*})\leftarrow \tilde{V}_{tt^*}(\tilde{S}_{tt^*})+\frac{\Delta- \tilde{V}_{tt^*}(\tilde{S}_{tt^*})}{N(\tilde{S}_{tt^*})}$
\STATE \hspace{1cm} $t'\leftarrow t^*$
\STATE \hspace{0.5cm}\textbf{end while}
\end{algorithmic}
\end{algorithm}

The choice of the parameters of MCTS is problem dependent; in this work, we choose the following parameters for the four steps of MCTS to find the optimal route of the truck.

\begin{itemize}
\item[1. ]\textbf{Selection:} In sampled MCTS, there are two selection strategies which are applied based on the domain of selection; one is  for the decision space while the other is for the exogenous event space. Starting from the root node, selection chooses  a decision based on previous gained information while controlling a  balance between exploration and exploitation. The most popular method used in the computer science literature is Upper Confidence Bounding applied to Trees (UCT)~(\citet{BPW12,KS06}). However, upper confidence bounds are used for maximization problems. In this work, the aim is to minimize the objective function which is equivalent to maximizing its negative value. UCT builds on an extensive literature in computer science on upper confidence bounding (UCB) policies for multiarmed bandit problems~(\citet{ABF02}). UCT selects the decision that maximizes the following equation:
    \begin{equation}
    \tilde{x}_{tt'}^*=\arg\max_{\tilde{x}_{tt'}\in \tilde{\mathcal{X}}_{tt'}^e(\tilde{S}_{tt'})}\left(-\left(\tilde{C}(\tilde{S}_{tt'},\tilde{x}_{tt'})+\tilde{V}^x_{tt'}(\tilde{S}^x_{tt'})\right)+\alpha\sqrt{\frac{\ln N(\tilde{S}_{tt'})}{N(\tilde{S}_{tt'},\tilde{x}_{tt'})}}\right),
    \end{equation}
    where $\alpha$ is a tunable parameter that balances exploration and exploitation. The choice of decision $\tilde{x}_{tt'}^*$ that maximizes the UCT equation depends on a weighted average of two terms; the first term of the UCT equation represents the average value of the state-decision after $N(\tilde{S}_{tt'},\tilde{x}_{tt'})$ iterations. So, the higher the average value of the state-decision, the more it contributes to exploiting the decision further since its reward is high. The second term gives a higher weight to the decision that has been less explored since its value decreases as $N(\tilde{S}_{tt'},\tilde{x}_{tt'})$ increases which contributes to exploring the decisions with lower number of visits.

 \indent Upon choosing $\tilde{x}_{tt'}^*$, the state of the network becomes $\tilde{S}^{x}_{tt'}$ after which we sample an exogenous realization $\tilde{W}_{t,t'+1}$ from the set of explored exogenous events $\tilde{\Omega}^e_{t,t'+1}(\tilde{S}_{tt'})$ for state $\tilde{S}_{tt'}$.
     In the developed model, the exogenous events may have very different probability density functions where some events can lie on the tail of the probability density function. Thus, in order to avoid too many iterations to catch the rare events, we propose importance sampling to choose $\tilde{W}_{tt'}^*$ from the set of available samples. Importance sampling yields the same expected value of the outcome of a random variable with a much lower number of iterations compared to sampling using the random variable's initial probability density function. Assume that $\tilde{\omega}$ represents an outcome of the exogenous random variable $\tilde{W}_{tt'}$. Since there are only a few outcomes for $\tilde{W}_{tt'}$, let $p(\tilde{W}_{tt'}=\tilde{\omega})$ be the probability mass function for outcome $\tilde{\omega}$ and $E_{p}[\tilde{W}_{tt'}]$ be its expected value. Also, define a new probability mass function $g(\tilde{W}_{tt'}'=\tilde{\omega})$ which is designed to balance the selection of all outcome events. According to importance sampling, $E_{p}[\tilde{W}_{tt'}]\approx E_{g}[p(\tilde{W}_{tt'}')/g(\tilde{W}_{tt'}')\tilde{W}_{tt'}']$ for a much lower number of realizations of $\tilde{W}_{tt'}'$. For simplicity, let $\tilde{\omega}_{tt'}$ be an abbreviation of the event $\tilde{W}_{tt'}=\tilde{\omega}$. In this work, we choose $g(\tilde{W}_{tt'})$ to have a uniform distribution for all random events that can take place from state $\tilde{S}_{t,t'-1}^{x}$. Then, one of the outcomes $\tilde{\omega}$ is chosen according to $g(\tilde{W}_{tt'})$ and later in the backpropagation step,  the value function of $\tilde{S}^{M,x}(\tilde{S}_{t,t'-1}^{x},\tilde{\omega}_{tt'}^*)$ is weighted by $p(\tilde{\omega}_{tt'}^*)/g(\tilde{\omega}_{tt'}^*)$ in order to maintain the same expected value of the exogenous events.

\item[2. ]\textbf{Expansion:} This is the process of adding a child node  to the tree to expand it. Upon visiting a state, one can either expand an unexplored state (via a decision or exogenous event) or exploit existing states. For example, one can set a threshold, $d^{thr}$, for the number of decisions and another threshold, $e^{thr}$, for the number of exogenous events to be expanded first before starting the exploitation process.   If a state has several unexplored decisions and exogenous events represented by the sets $\tilde{\mathcal{X}}_{tt'}^u(\tilde{S}_{tt'})$ and $\tilde{\Omega}_{tt'}^u(\tilde{S}_{tt'}^{x})$, respectively,  then an unexplored decision is chosen optimistically by using a two-stage lookahead model. That is, for each unexplored decision $x_{tt'}\in \tilde{\mathcal{X}}_{tt'}^u(\tilde{S}_{tt'})$, we generate a random exogenous event $\tilde{\omega}\in \tilde{\Omega}_{t,t'+1}^u(\tilde{S}_{tt'}^{x})$ and evaluate the value of the obtained state $\tilde{S}^{M,x}(\tilde{S}_{tt'}^{x},\tilde{\omega}_{t,t'+1})$ by calling the simulation policy. Then, the decision and its corresponding exogenous event that corresponded to the highest obtained value are chosen to be expanded in the tree. Finally, if the selection phase reaches a  state, say  $\tilde{S}_{tt'}^{x}$, that is already part of the tree but for which the threshold value of the number of exogenous events to be explores is not met, then a random exogenous sample $\tilde{\omega}\in \tilde{\Omega}_{t,t'+1}^u(\tilde{S}_{tt'}^{x})$ is created resulting in state  $\tilde{S}^{M,x}(\tilde{S}_{tt'}^{x},\tilde{\omega}_{t,t'+1})$ and evaluated as discussed above.

\item[3. ]\textbf{Simulation:} The simulation policy is a heuristic policy to provide an initial estimate of the value of the state that has just been added to the tree. Adding a node to the search tree at time $t'$, results in state $\tilde{S}_{tt'}$ and updates the probability space $\tilde{\Omega}_{tt'}(\tilde{S}_{tt'})$. The simulation policy selects states from the newly added state until the end of the simulation; this is a roll-out simulation starting from the expanded state. To get an efficient simulation policy, we propose one that it is not too stochastic (for efficiency of evaluation) nor too deterministic (because this will bias the search tree)~(\citet{BPW12}).
    The proposed simulation policy is based on a lookahead model in which we evaluate the value of the newly added state at time $t'$ by generating a sample path $\tilde{\omega}\in\tilde{\Omega}_{tt'}(\tilde{S}_{tt'})$. The sample path determines the set of power lines  that have faulted along with the fault types, and required travel and repair times for each arc in the graph.
    We elaborate on the simulation policy in Appendix~\ref{sec:sim_policy} where we formulate the problem as a sequential optimization problem; the resulting problem is a mixed integer non-linear program (MINLP) where the nonlinearity arises from the radial structure of the grid. In general, MINLP problems are difficult to solve for large network sizes but Appendix~\ref{sec:sim_policy} shows that the problem reduces to  a travelling salesman problem (TSP) since the objective is to find the tour of the truck that visits each location with a fault exactly once to repair it in order to minimize the customer outage-minutes. So, Appendix~\ref{sec:sim_policy} shows that the optimal route of the truck given a sample path can be optimally solved via dynamic programming for a small number of generated faults whereas a heuristic TSP solution becomes necessary if the number of generated faults is large.

\item[4. ]\textbf{Backpropagation:} At the end of the simulation, a value, $\tilde{V}_{tt'}(\tilde{S}_{tt'})$, of the newly created state, $\tilde{S}_{tt'}$, is obtained. Starting from the last added node in the tree, its simulated value is backpropagated through all ancestors of state $\tilde{S}_{tt'}$  until the root state to update their statistics. The number of visit counts of all  ancestor states of state $\tilde{S}_{tt'}$ are increased by one and their values are modified according to a chosen criteria where we choose the average value of all rollouts through a state. While backpropagating, assume that we have $\tilde{S}_{tt'}=S^{M,x}(\tilde{S}_{tt^*}^x,{\tilde{\omega}_{t,t^*+1}})$, then the value of the ancestor post-decision state, say $\tilde{S}_{tt^*}^x$, should be updated as $\tilde{V}^x_{tt^*}(\tilde{S}^x_{tt^*})\leftarrow \frac{1}{\sum_{\tilde{W}_{t,t^*+1}\in\tilde{\Omega}_{t,t^*+1}^e(\tilde{S}_{tt^*}^{x})} p(\tilde{W}_{t,t^*+1})} \cdot E_{g}[p(\tilde{W}_{t,t^*+1})/g(\tilde{W}_{t,t^*+1})\tilde{V}_{tt^*}(S^{M,x}(\tilde{S}_{tt^*}^x,{\tilde{W}_{t,t^*+1}}))]$ where the expectation is over all the explored exogenous events from post-decision state $\tilde{S}_{tt^*}^x$. Also, the value function of the ancestor pre-decision state, $\tilde{S}_{tt^*}$ should be updated with a value of $\Delta \leftarrow \tilde{C}(\tilde{S}_{tt^*},{\tilde{x}_{tt^*}}) + \tilde{V}^x_{tt^*}(\tilde{S}^x_{tt^*})$ which accounts for the link cost and the updated post-decision state value function so that we get
  $\tilde{V}_{tt^*}(\tilde{S}_{tt^*})\leftarrow \tilde{V}_{tt^*}(\tilde{S}_{tt^*})+\frac{\Delta- \tilde{V}_{tt^*}(\tilde{S}_{tt^*})}{N(\tilde{S}_{tt^*})}$ (note that we assume that the weight of all decisions from the same state is equal).
\end{itemize}

\section{Performance Results}\label{sec:results}
To assess the performance of the proposed approaches, the simulated power grid is constructed using real data provided by PSE\&G which describes the structure of circuits in their electrical distribution network. The data corresponds to the northeastern portion of
PSE\&G's power grid in New Jersey and it is formed of $319$ circuits. There are an average of $41$ protective devices and $724$ power lines per circuit. The data identifies the type and location of each component in the circuits such as substations, protective devices, power lines and transformers.

The simulator is programmed to generate storms that pass across the grid generating power line faults causing total or partial circuit power outages~(\citet{ABP14}). The obtained outages trigger some of the affected customers to call to report the outage. When a power line faults due to storm damage, the simulator finds the nearest upstream protective device, opens it (shutting off power), and then identifies all the customers who subsequently lose power. For each customer experiencing a power outage, a Bernoulli random variable is generated, with a probability of success which equals to the calling probability across the customer's segment,  to determine whether the customer will call or not. Thus, the higher the calling probability is, the higher the number of trouble calls will~be.

In order to reconstruct the grid due to storm damage, a truck uses  a roadway defined by the graph $G(\mathcal{V},\mathcal{E})$ where $\mathcal{E}$ is the set containing the arcs/roads between two consecutive nodes of the graph. 
We represent the network by aggregating lines and poles that are protected by the same protective device into a single segment.  The idea is that if a truck visits one pole or line of the circuit, it will be able to quickly see the status of nearby lines and poles, which we represent as a segment.  We then assume that we are guiding trucks from one segment to another (where the truck has to find the best path over the road network).  When a truck visits a segment, it is assumed to learn (and fix, if necessary) all outages anywhere in a segment.
Basically, the network  $G(\mathcal{V},\mathcal{E})$ is formed of the poles that carries the protective devices where $\mathcal{E}$ corresponds to the connection matrix in the form of the minimum distance between any two nodes in the network according to the real roadway.

As the storm passes across the grid, each power line $i$ along its way is associated with a prior probability of power line fault as explained in detail in~(\citet{ABP14}). In the simulator, the prior probability of fault of the power lines are generated based on several parameters such as the severity of the storm, its diameter and the distance of the power line from its center (refer to~(\citet{ABP14}) for more details). Since this paper addresses a single truck, we tuned the priors to create storms that generated one to as many as tens of outages. In this case, the segment fault probabilities range between  $0$ and $0.765$ where the segments that faulted have posteriors ranging between $0.032$ and $0.765$.

\subsection{Lookahead Policy}
After collecting the priors for power line faults, the obtained faults and the customers that called, the simulator executes the proposed power line fault probability model presented in (\ref{eq:eq3}) upon which the utility truck is routed to restore the grid.
 According to the simulations, the average number of segments affected by the storm path is $1558$; however, the probability model sets the posterior probability of a major number of them to a very low number, e.g., below $0.01$. Based on the statistics of $1000$ networks, the minimum posterior fault of a segment that faulted is $0.032$, so we choose the segments with posterior probability of faults greater or equal to $0.01$
as candidate segments that have faulted.

\begin{figure}[t!]
\centering
\includegraphics[scale=0.58]{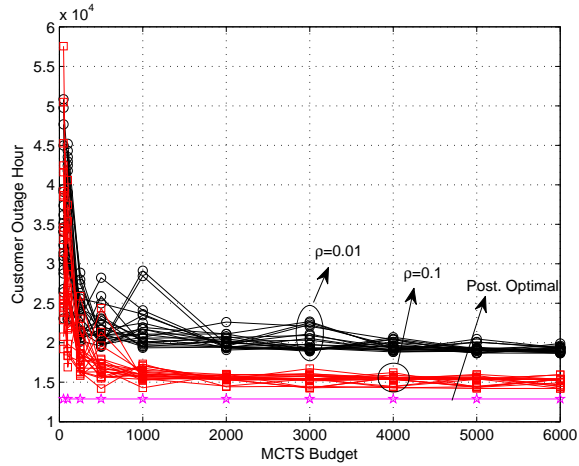}
\caption{\label{fig:Simu_COH} Customer outage-hours vs. MCTS budget for one network; twenty simulations for various calling probabilities where each simulation corresponds to the route of the truck starting from the depot till the stopping condition is met.}
\end{figure}

Figures~\ref{fig:Simu_COH} shows the customer outage-hours vs. MCTS budget for twenty simulations for one chosen network which has $5$ faults. By MCTS budget, it is meant the number of iterations, $n^{thr}$,  executed to return the decision of where to route the truck next. A simulation represents the route of the truck from the depot until the stopping condition of truck routing presented in Algorithm~\ref{alg:adp} is met. There are two stopping conditions; one to stop MCTS search when the number of iterations performed, to return a final decision, reaches  $n^{thr}$. The other is for stopping truck routing when the posterior probabilities of faults of all nodes in the network drops below a threshold $\epsilon^{thr}$ which is set to $0.01$ since the posterior probabilities of the segments that faulted are much higher.

In MCTS, there are several steps that depend on random outcomes such as which exogenous event to choose or the generated sample path to return a value function for a newly generated state; this results in different solutions for different calls of MCTS especially for a low budget; however, as the MCTS budget increases, the averaging over the random events becomes more accurate and MCTS returns solutions with nearly similar objective values. Figure~\ref{fig:Simu_COH} also shows that the rate of convergence to a similar objective value depends on the calling probability; that is, as the calling probability increases, then more information is revealed in the network which reduces the uncertainty and consequently MCTS returns a good solution with a lower budget. According to Figure~\ref{fig:Simu_COH}, a low MCTS budget is not enough to decide on the route of the truck but also a very high MCTS budget does not provide any improvement at the cost of higher computational complexity. The simulations are also compared to the posterior optimal solution which corresponds to the optimal solution after revealing the locations of faults in the network  obtained using the dynamic program presented in Algorithm~\ref{alg:dptsp}.

Also, one of the important tunable parameters in the MCTS algorithm is the parameter that balances the weight between exploitation and exploration (referred to as $\alpha$ in Algorithm~\ref{alg:MCTS}). For the same  network considered above, a brute force search over the best value of $\alpha$ is conducted ranging from $0.1$ till $7$ with increments of $0.1$. Starting from $0.1$, the average objective value decreases reaching its minimum at $\alpha=2.2$ and then it increases again. So, for the rest of the simulations, we set $\alpha = 2.2$.

\begin{figure}[t!]
\centering
\includegraphics[scale=0.58]{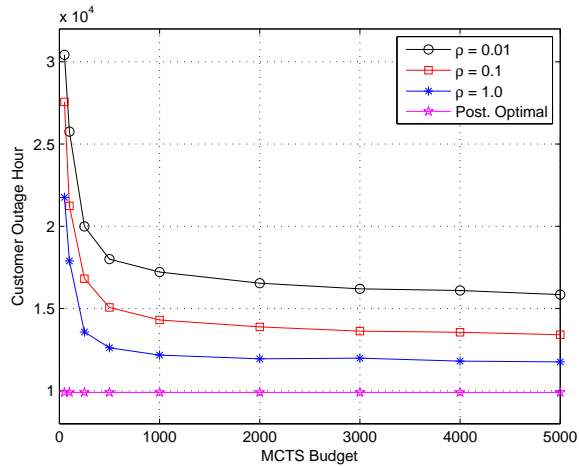}
\caption{\label{fig:Av_COH}Average customer outage-hours vs. MCTS budget for ten networks.}
\end{figure}

Figure~\ref{fig:Av_COH} shows the average customer outage-hours vs. MCTS budget of ten networks for various calling probabilities. The number of faults in the networks ranges from $4$ to $12$ with an average of $6.09$ faults. First, as the  calling probability increases, the lookahead policy provides a solution that is closer to the posterior optimal since more information is provided. Also, the required MCTS budget decreases as the calling probability increases; that is, for a calling probability of $0.01$ and $1.0$, around $4000$ and $1000$ MCTS iterations are required, respectively, to converge to a good solution. The lookahead policy provides a solution that is $18.7$\%, $35.3$\% and  $58.5$\% higher than the posterior optimal for a calling probability of $1.0$, $0.1$ and $0.01$ respectively. Second, it is revealed that a high gain is obtained when the calling probability increases from $0.01$ to $0.1$ since the information provided can be from different locations which help in detecting the location of outages. As the calling probability rises above $0.1$, the benefits are more modest than the increase from $0.01$ to $0.1$ since we only need one customer out of a group across the same segment to make the phone call.

\begin{figure}[t!]
\centering
\includegraphics[scale=0.58]{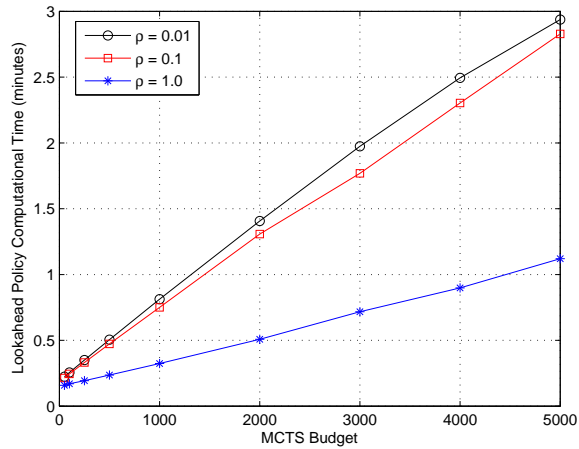}
\caption{\label{fig:Av_CT}Average computational time of the lookahead policy to decide on the truck's next hop vs. MCTS budget.}\vspace{-0.4cm}
\end{figure}

Figure~\ref{fig:Av_CT} shows the average computational time of the lookahead policy to decide on the next truck route. Obviously, the computational time increases as the MCTS budget increases; however,  no additional gain in terms of the objective value is obtained.  It is also shown that the computational time decreases as the calling probability increases since more information is provided about the locations that could have faulted  which reduces the network size after considering all nodes with posterior probabilities above a certain threshold. However, in the worst case, for a calling probability of $0.01$ and a good computational budget of MCTS, which is shown to be $4000$ in Figure~\ref{fig:Av_COH}, the computational time of the lookahead policy is $2.5$ minutes which is suitable for online problems.

\begin{figure}[t!]
\centering
\includegraphics[scale=0.58]{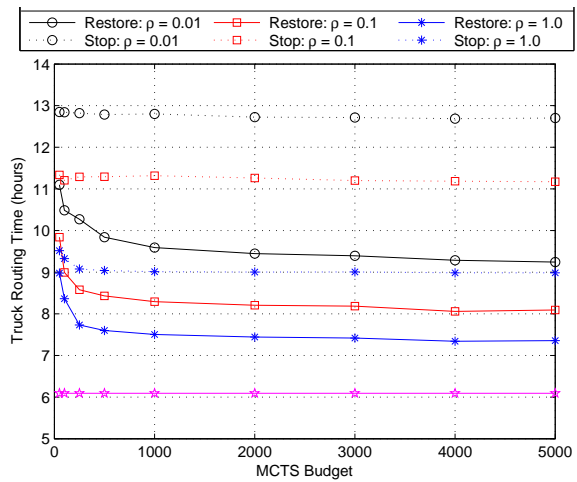}
\caption{\label{fig:AV_TRT} Average truck routing time to restore the grid and to stop routing vs. MCTS budget; the truck is stopped when the stopping condition of the lookahead policy is met.}\vspace{-0.4cm}
\end{figure}

Figure~\ref{fig:AV_TRT} shows the time required to restore the grid and the time required to stop truck routing compared to the posterior optimal solution.  Even when the grid is completely restored,  the utility center cannot  detect that unless the distribution system is fully equipped with SCADA which is not the case. In this work, the utility center relies on the probability model to decide when the grid is restored and consequently, stop the truck routing process.  The posterior optimal solution indicates that  $6$ hours are required to restore the grid if full information is revealed. For a calling probability of $1.0$, the lookahead policy can restore the grid in $7.2$ hours on average but, the truck is stopped after $11$ hours. For a lower calling probability, both the restore and stop times increase for the same reasons mentioned previously but still have a good performance with respect to the posterior optimal~solution.

%

\subsection{Industrial Heuristics}
We also compare the performance of the proposed  lookahead policy to industrial heuristics which typically rely on escalation algorithms. An escalation algorithm back-traces from each trouble call location to find the first common point for all the calls. It is obvious that escalation algorithms are good to locate a single fault which is assumed to be an upstream fault that triggers all downstream calls. But, often this is not the case as there can be more than one fault triggering the calls.  Practically, escalation is performed at the control center along with other intelligence techniques as explained in the literature review in  Section~\ref{sec:lit_rev}. However, in this work, only the trouble calls and the grid structure are used as the sources of information. Thus, the proposed escalation algorithm gives priority to visit the common locations that would trigger all calls but after that searches for other downstream faults as explained in Algorithm~\ref{alg:esc}.

\begin{algorithm}[t!] \footnotesize
\caption{Escalation Algorithm for Grid Restoration}
\label{alg:esc}
\begin{algorithmic}[t!]
\STATE\textbf{Step 1.} \textbf{For }each circuit \textbf{do}
\STATE\hspace{1.1cm}\textbf{Step 1a.} Collect all calls and back trace to find the first node that is common to all calls say node $x$.
\STATE\hspace{1.1cm}\textbf{Step 1b.} Send the truck to node $x$ and then back trace to the substation to cover all upstream faults \STATE\hspace{1.5cm}(when a truck visits a node, it fixes an existing fault and this applies to all steps of the algorithm).
\STATE\hspace{1.1cm} \textbf{Step 1c.} From node $x$, perform down tracing to reach the first segment from which a call was initiated \STATE\hspace{2.5cm}  and place it in set~$\mathcal{D}$.
\STATE \textbf{Step 2.} \textbf{For} each segment in $\mathcal{D}$ \textbf{do}
\STATE\hspace{1.1cm} \textbf{Step 2a.} Perform down tracing to cover all nodes that called.
\end{algorithmic}\vspace{-0.0cm}
\end{algorithm}

\begin{table}[t]
\vspace{-0.0cm}
\caption{Average Statistics for Escalation Algorithm}\label{tab:esc}
\begin{center}\vspace{-0cm}
\begin{tabular}{|l|c|c|c|}
\hline
 & $\rho = 0.01$ & $\rho = 0.1$ & $\rho = 1.0$   \\
\hline
Customer outage-hours & $ 2.34*10^4$& $2.21*10^4$ & $2.17*10^4$\\
\hline
Number of unrepaired faults& $1$ &$0.42$ &$0$ \\
\hline
Number of customers in outage & $14$ & $ 5.9$ & $0$ \\
\hline
Time to restore (hours)& $21.83$ &  $20.76$ & $20.34$ \\
\hline
Time to stop (hours) & $48$ & $48$ & $26.46$ \\
\hline
\end{tabular}
\end{center}\vspace{-0.5cm}
\end{table}

Table~\ref{tab:esc} shows the average statistics of the same ten networks used earlier. In the escalation algorithm, searching for a fault is just triggered by the customer calls because there is no clue for the utility center to predict the locations of faults except by visiting the locations that  trigger the calls if faulted. So, if there is a segment with a low number of customers where no one called to report an outage then, there is no way to tell that a fault exists unlike the proposed probability model which would recognize that there is a nonzero probability of an outage on these segments due to the storm path.

In Table~\ref{tab:esc}, we see that for a calling probability of $0.01$, on average $1$ fault could not be identified since a total of $14$ customers are only affected and non of them called at such a low calling probability. So, the utility center will stop routing the truck assuming that it has recovered the grid, because it visited every location from which a call was received. The average number of faults is $6.09$ with an average of $1$ unrepaired fault. So, the time it to took the truck to restore an average of $5.09$ faults that triggered calls is around $21.83$ hours and the customer outage-hours is $2.34*10^4$. We set the maximum truck routing time to stop to $48$ hours unless it restored all the grid in a smaller amount of time. The escalation algorithm was fast in locating an upstream fault whenever there is one that triggered all downstream calls but it took it a long time to locate downstream faults because there is no clue from the common point and on except by following the locations that could trigger the calls. Whereas for the same calling probability, the lookahead policy is able to repair all faults in $9.15$ hours with a customer outage-hours of $1.58*10^4$ which means that is much better because it was able to reconstruct all faults with a smaller objective function and truck routing time.

The lookahead policy outperforms the escalation algorithm by orders of magnitude as the calling probability increases; for example, for a calling probability of $1.0$, the escalation algorithm is able to identify all fault locations because basically every customer that lost power is calling. But, the restore and stop times are $20.34$ and $26.46$ hours, respectively, whereas the customer outage-hours is $2.17*10^4$. More calls can help identify upstream faults extremely fast but there is no strategy for identifying the downstream faults because everyone is calling due to the upstream faults. So, increasing the calling probability is not as beneficial for the escalation policy as it is for the lookahead policy which makes a much clever use of the provided information by decreasing the restore and repair times to $7.2$ and $9$ hours with an objective of $1.17 * 10^4$. 

Based on the results, we see that the current rules followed in industry can work well for locating upstream faults but it does not provide accurate guidance for finding cascaded downstream faults. Moreover, if no calls are received from a location with an outage, then the escalation algorithm cannot find it. Not only that, when increasing the calling probability, the escalation does not make good use of the additional information except for locating upstream faults quickly.


\section{Conclusion}\label{sec:conclusion}
In this work, a lookahead policy is proposed to route the  utility truck across the power grid to restore it as quickly as possible after a storm. In addition to the trouble calls, the utility truck was used as a mechanism for collecting additional information to update the belief about faults in the grid. Performance results have shown that even with a small percentage of customers calling to report an outage,  the lookahead policy can restore the power grid efficiently. It is also shown that the lookahead policy outperforms current techniques used in industry that are based on escalation to locate the faults and in making use of the additional provided information.
\setlength{\bibsep}{1.0pt}
\bibliographystyle{elsarticle-harv}
\bibliography{references}

\begin{APPENDICES}
\section{Lookahead Simulation Policy}\label{sec:sim_policy}
In the lookahead simulation policy, we evaluate the value of the newly added node at time $t'$ by using a lookahead model in which a sample path $\tilde{\omega}\in\tilde{\Omega}_{t,t'+1}(\tilde{S}_{tt'})$  is generated. The sample path determines the set of power lines  that have faulted along with the fault types, required repair and travel times for each arc in the graph.  For example, using the probability of fault of each power line in the system, we generate a Bernoulli random variable with a probability of success equals to the probability of fault. In this case, each power line in the power system has a posterior probability of fault equal to either $1$ or $0$. Thus, we define the following indicator function:
\begin{eqnarray}
\mathds{1}_{\tilde{L}^u_{t'i}}=\left\{
\begin{array}{c l}
    1 & \mbox{,~if~} \tilde{L}^u_{tt'i}(\tilde{\omega})=1,\\
    0 & \mbox{, otherwise},
\end{array}\right.
\end{eqnarray}
where $\mathds{1}_{\tilde{L}^u_{t'i}}=1$ if  power line $i$ on circuit $u$ is in fault at time $t'$ and it is equal to zero, otherwise.

In the simulation policy, we use the time index $t''$ for each state, decision and random variable. Recall, at time~$t$, we are in the base model where we fix the set of calls and call MCTS to find the truck's next hop. For each node included in the MCTS tree, we index it with time $tt'$.
In the fourth step of MCTS, to evaluate the value of a newly added state $\tilde{S}_{tt'}$, we use another lookahead model to generate a sample path $\tilde{\omega}\in\tilde{\Omega}_{t,t'+1}(\tilde{S}_{tt'})$ at $t'$. While, routing the truck according to the generated path, we index the nodes  by $t''$; for example, $\tilde{S}_{t''}$ indicates the state at time $t''$ in the simulation step which is used to evaluate the value of the generated state in MCTS at time $t'$.

In the simulation step, if at time $t''$, the truck visits a location that was identified as a location with fault at time $t'$, its indicator function is set to $0$ from time $t''+1$ and on. Thus, for each power line that faulted  the following relation holds:
\begin{eqnarray}
\mathds{1}_{\tilde{L}^u_{t''j}}=1-\sum_{i\in\mathcal{V}}\sum_{\hat{t}=t'}^{t''-1}\tilde{x}_{\hat{t}ij}, \mbox{~if~} \tilde{L}^u_{tt'j}(\tilde{\omega} )=1.\label{eq:const5ex}
\end{eqnarray}

Whereas, if power line $i$ on circuit $u$ did not fault in the considered scenario then, $\mathds{1}_{\tilde{L}^u_{t''i}}=0$ , $\forall t''\geq t'$.
Given, a sample path $\tilde{\omega}$, the aim is to find the optimal truck's route that minimizes the customer outage-minutes. Let $\tilde{C}_{t''j}$ be a random variable representing the number of customers that regain power by visiting node $j$ at time $t''$ according to scenario $\tilde{\omega}$. Then, the value of the newly added node is obtained by solving the  following optimization problem:

\begin{eqnarray}
&&\hspace{-1.0cm} \min_{\tilde{x}_{t''}} \sum_{\hat{t}=t'}^{t'+H}\left(N-\sum_{t''=t'}^{\hat{t}}\sum_{u\in\mathcal{U}}\sum_{j=1}^N\tilde{C}^u_{t''j}\right)\label{eq:objective}\\
&&\hspace{-1cm} \mbox{subject to}\nonumber\\
&&\hspace{-1cm}\tilde{C}^u_{t''j}=\sum_i \left(\prod_{k\in\mathcal{Q}^u_j\backslash j}\hspace{-0.2cm}1-\mathds{1}_{\tilde{L}^u_{t''k}}\right)\mathds{1}_{\tilde{L}^u_{t''j}}\left(\sum_{k\in \mathcal{S}^u_j} n^u_k + \sum_{s\in \mathcal{W}^u_{j}}\left(\prod_{w=\min\{\mathcal{W}^u_{j}\}}^s \prod_{k\in w}1-\mathds{1}_{\tilde{L}^u_{t''k}}\right)\sum_{k\in s}n^u_k\right)\tilde{x}_{t''ij}, \nonumber\\ &&\hspace{12.5cm}\forall j\in\mathcal{V}, \forall t''\\ \label{eq:const1} 
&&\hspace{-1cm} \mathds{1}_{\tilde{L}^u_{t''j}}=1-\sum_i\sum_{\hat{t}=t'}^{t''-1}\tilde{x}_{\hat{t}ij},\mbox{~such~that~} \tilde{L}^u_{tt'j}(\tilde{\omega})=1 ,\forall j \in \mathcal{I}^u,\forall u\in\mathcal{U}, \forall t''\label{eq:const5}\\
&&\hspace{-1cm}\tilde{\Delta}_{t''ij}\geq T_{ij}(\tilde{\omega})\tilde{x}_{t''ij} + \sum_u R^u_{j}(\tilde{\omega}) \left(\tilde{x}_{t''ij} -\sum_{\hat{t}=t'}^{t''-1}\tilde{x}_{\hat{t}ji} \right),\forall (i,j)\in \mathcal{E},\forall t''\label{eq:const2}\\
&&\hspace{-1cm}\tilde{\xi}_{t''j}\geq \tilde{\xi}_{t''-1i}+ \sum_i \tilde{\Delta}_{t''ij} ,\forall (i,j)\in\mathcal{E}, \forall t''\label{eq:const3}\\
&&\hspace{-1cm}\tilde{\xi}_{t''j}\leq t''\sum_i\tilde{x}_{t''ij}+\zeta\left(1-\sum_i\tilde{x}_{t''ij}\right),\forall j\in\mathcal{V}, \forall t''\label{eq:const4}\\
&&\hspace{-1cm}\sum_{t''=t'}^{t'+H}\tilde{x}_{t''ij}\leq 1,\forall (i,j)\in \mathcal{E}\label{eq:const7}\\
&&\hspace{-1cm} \sum_k \tilde{x}_{(t''+T_{jk}(\tilde{\omega}))jk}+\sum_k \tilde{x}_{(t''+T_{jk}(\tilde{\omega})+\sum_uR^u_{k}(\tilde{\omega}))jk}\leq \sum_{i}\tilde{x}_{t''ij}\leq 1,\forall j\in\mathcal{V},\forall t''\label{eq:const8}\\
&&\hspace{-1cm} \tilde{C}^u_{t''j}\geq 0, \tilde{\xi}_{t''j}\geq0, \tilde{\Delta}_{t''ij}\geq0, \tilde{x}_{t''ij}\in\{0,1\} \label{eq:const11}
\end{eqnarray}

This problem is a mixed integer non-linear programming (MINLP) problem. The objective~(\ref{eq:objective}) minimizes the customer outage-minute which is equivalent to maximizing the number of customers with restored power (also referred to as served customers) up to time $t$ represented by the inner sum in the objective.
Constraint (\ref{eq:const1}) determines the  number of served customers when the truck goes from its current  location, say node $i$, to node $j$ at time $t''$. The number of served customers depends on whether there is a fault upstream to node $j$ or on its segment, i.e., in set $\mathcal{Q}_j^u$. Note that, node $j$ will be favored to be visited if there is a fault across power line $j$, i.e., if $\mathds{1}_{\tilde{L}^u_{t''j}}=1$.
Depending on the structure of the power grid, if a location faults, it causes outage to the customers attached to its segment and all downstream segments. Thus, if a fault is fixed, then all these customers will be affected. But, this also depends on whether there is a fault on any downstream location as shown in~(\ref{eq:const1}). Constraint~(\ref{eq:const1}) also reveals that the number of customers by visiting power line $j$  is positive if $\mathds{1}_{\tilde{L}^u_{t''j}}=1$; however, after visiting this location, say at time $t^*$, $\mathds{1}_{\tilde{L}^u_{t''j}}=0, \mbox{~for~} t''> t^*$. Thus, if the truck will come across the same location for the second time, then the gain will be $0$ which favors the truck not to visit the same location more than once unless there is no other route for it. Constraint (\ref{eq:const5}) is the same as (\ref{eq:const5ex}) and it has been explained in details above.

Constraint (\ref{eq:const2}) determines the required time to traverse arc $(i,j)\in \mathcal{E}$. If there is no power line to be investigated at node $j$, then $\mathds{1}_{\tilde{L}^u_{t''j}}=0$ which means that the required traversal time is equal to the  travel time which depends on the traffic conditions only. However, if there is a power line across arc $(i,j)$, then one of two rules apply; if there is a fault on power line $j$ according to sample path $\tilde{\omega}$,  then the required traversal time accounts for the travel and repair times for power line $j$. However, if arc $(i,j)$ is traversed for the second time at time $t$ then, the traversal time is just equal to the travel time since the fault was repaired when the arc was traversed for the first time.

Constraints (\ref{eq:const3})-(\ref{eq:const4}) guarantee that the truck is at node $j$ at time $t''$, only if $\tilde{\xi}_{t''j}=t''$ which sets $\tilde{x}_{t''ij}=1$. If $\tilde{x}_{t''ij}=1$, then $\tilde{\xi}_{t''j}=t''$, otherwise $\tilde{\xi}_{t''j}$ is less than a large positive number $\zeta$ as shown in (\ref{eq:const4}) but larger than the time where the truck was lastly as indicated by~(\ref{eq:const3}). But, since the objective is maximizing the number of served customers over time, the optimization problem will set the time to the least possible value that satisfies all constraints.  In (\ref{eq:const3}), $\tilde{\xi}_{t''j}$ can be equal to $t''$ only if it satisfies the required traverse times; the required time to reach node $j$ depends on the elapsed time to reach its direct predecessor, say node $i$, in addition to the required time to traverse node $j$ from node $i$, i.e., $\tilde{\Delta}_{t''ij}$.

Note that, since the objective minimizes the customer outage-minutes, then the optimization problem keeps on routing the truck to cover all power lines that faulted, as favored by (\ref{eq:const1}), until all faults are fixed. Constraint (\ref{eq:const7}) guarantees that all arcs in the graph can be visited once in one direction and consequently, at most twice (forward and backward) which is a sufficient condition to have an Eulerian path where each power line and node with positive fault probability can be visited once. Constraint~(\ref{eq:const8}) indicates that the truck can go from node $j$ to node $k$ at time $t''+\tilde{\Delta}_{t''jk}$ only if it was at node $j$ at time $t''$ and only if the necessary traverse time, $\tilde{\Delta}_{t''jk}$, has elapsed  which depends on (\ref{eq:const2}). Moreover, this constraint removes the sub-tours in the network since the truck must have visited a node before it can travel from~it. Finally, constraint~(\ref{eq:const11}) shows that all variables are positive except $\tilde{x}_{t''ij}$ which is binary.

The formulated  optimization problem is very complex mainly because it is an MINLP problem which combines the complexity of non-linear programming and integer programming both of which lie in the class of NP-hard problems.  Thus, achieving the optimal global solution is most probably never attainable for large network sizes. While there has been a tremendous achievements in solving integer programming problems given that their continuous relaxation is convex, solving non-linear optimization problems is still a non mature area that gets stuck at local optimums.

The only constraint that cannot be linearized is (\ref{eq:const1}); it can be seen that the order of non-linearity depends on the number of faults upstream and downstream of a node which is scenario dependent. Thus,  the radial structure of the grid   is the main complicating factor in the optimization problem. 

Though the problem is non-linear, the optimal solution can be attained by using dynamic programming. For each sample path, we can transform the problem into a complete graph with nodes $\mathcal{V}^f$ which contains all the power lines that have faulted and the location of the truck indexed with $0$. The connection cost between the nodes of $\mathcal{V}^f$ are calculated by summing the shortest travel time between the nodes according to $\tilde{\omega}$. Let $S$ be the set of the nodes visited by the truck and $f(S)$ a function returning the number of customers still in outage after visiting the nodes of $S$. The aim of the problem is to find the optimal sequence of the truck route that visits each node exactly once (to repair it) in order to minimize the customer outage-minutes. This problem is equivalent to a travelling salesman problem (TSP) which is NP-Complete; however the solution can be obtained optimally using dynamic programming with complexity $O(n^22^n)$ where $n$ is the number of nodes in the TSP graph which corresponds to the number of generated faults. However, since the number of generated faults is relatively small (less than 20 faults), obtaining the optimal solution using dynamic programming is computationally  feasible.

Let $C^f(S,i)$ be the customer outage-minutes of going from vertex $0$ through the nodes of $S$ ending at node~$i$. Then, the recurrence relation of the dynamic program by going from node~$i$ to node~$j$ can be defined as
\begin{eqnarray}
C^f(S,j) =  C^f(S-\{j\},i)+f(\{S-\{j\}\})\cdot(T_{ij}+\sum_uR^u_j),
\end{eqnarray}
where the first term of the summation accounts for the customer outage-minutes up to node~$i$ whereas the second term accounts for the  cumulative customer outage-minutes by going from node~$i$ to node~$j$. The detailed steps of the dynamic program to obtain the value of the objective function are presented in Algorithm~\ref{alg:dptsp}.

\begin{algorithm}[t!] \footnotesize
\caption{Dynamic program for optimal customer outage-minutes of a given sample path}
\label{alg:dptsp}
\begin{algorithmic}[t!]
\STATE\textbf{Step 0.} \textbf{Initialization}
\STATE \hspace{1.2cm} \textbf{For }all $j\in\mathcal{V}^f, j\neq 0$ \textbf{do}
\STATE \hspace{1.4cm} $C^f(\{0,j\},j) = f(\{0,j\})\cdot(T_{0j}+\sum_uR^u_j)$
\STATE\textbf{Step 1.}~ \textbf{Compute customer outage-minutes for all subsets}
\STATE \hspace{1.2cm}\textbf{For} $s=3$ \textbf{to }$|\mathcal{V}^f|$
\STATE \hspace{1.2cm} \textbf{For }all subset of $\mathcal{V}^f$ of size $s$ \textbf{do}
\STATE \hspace{1.4cm} \textbf{For} all $j\in S, j\neq 0$
\STATE \hspace{1.6cm}   $C^f(S,j) = \min_{i\in S, i\neq j} C^f(S-\{j\},i)+f(\{S-\{j\}\})\cdot(T_{ij}+\sum_uR^u_j)$
\STATE \textbf{Step 2.} \textbf{Optimal solution}
\STATE \hspace{1.2cm}  $\min_{j\in\mathcal{V}^f}C^f(\mathcal{V}^f,j)$
\end{algorithmic}\vspace{-0.0cm}
\end{algorithm}

\end{APPENDICES}
\end{document}